\documentclass[final,3p]{elsarticle}
 \usepackage[latin1]{inputenc}
 \usepackage{graphics}
 \usepackage{graphicx}
 \usepackage{epsfig}
\usepackage{amssymb}
 \usepackage{amsthm}
 \usepackage{lineno}
 \usepackage{amsmath}
   \numberwithin{equation}{section}
\usepackage{mathrsfs}

\newtheorem{thm}{Theorem}[section]

\newtheorem{lem}[thm]{Lemma}
\newtheorem{prop}[thm]{Proposition}
\newtheorem{defn}[thm]{Definition}

\biboptions{sort&compress,square}
\allowdisplaybreaks

\begin{document}
\begin{frontmatter}
\author[rvt1]{Jian Wang}
\ead{wangj484@nenu.edu.cn}
\author[rvt2]{Yong Wang\corref{cor2}}
\ead{wangy581@nenu.edu.cn}
\cortext[cor2]{Corresponding author.}
\address[rvt1]{School of Science, Tianjin University of Technology and Education, Tianjin, 300222, P.R.China}
\address[rvt2]{School of Mathematics and Statistics, Northeast Normal University,
Changchun, 130024, P.R.China}

\title{The spectral torsion for the one form rescaled Dirac operator}
\begin{abstract}
The spectral torsion is defined by three vector fields and Dirac operators and the noncommutative residue.
 Motivated by the spectral torsion and the one form rescaled Dirac operator,
  we give some new spectral  torsion  which is the extension of spectral  torsion  for Dirac operators, and compute the spectral torsion  for the one form rescaled Dirac operator on even-dimensional spin manifolds without boundary.
\end{abstract}
\end{frontmatter}

\section{Introduction}
\label{1}

In Connes' program of noncommutative geometry \cite{Co1,Co2,Co3,Co4,Co5}, the role of geometrical
objects is played by spectral triples $(\mathcal{A},\mathcal{H},D)$,
where ${\mathcal{A}}$ is a noncommutative algebra with involution $\ast$,
acting in the Hilbert space ${\mathcal{H}}$ while $D$ is a self-adjoint
operator with compact resolvent and such that,
$
\lbrack D,a]\ \hbox{is bounded}\ \forall\,a\in{\mathcal{A}}\,.
$
 Similar to the commutative case and the canonical spectral triple $(C^{\infty}(M),\mathcal{L}^{2}(S),D)$, where $(M,g,S)$ is a closed spin manifold and $D$ is the Dirac operator acting on the spinor bundle $S$, the spectrum of the Dirac operator $D$ of a spectral triple $(\mathcal{A},\mathcal{H},D)$ encodes the geometrical information of the spectral triple. However, to gain access to this information, one
should first find a spectral formulation of the specific geometric notion, and then
extend it to the level of spectral triples.
In the general case,
for an arbitrary spectral triple   $(\mathcal{A},\mathcal{H},D)$
 and a self-adjoint element $h = h^* \in  \mathcal{A}$,  Connes and Moscovici \cite{Co7}
gave the properties of the `perturbed' triple
\begin{equation}
(\mathcal{A},\mathcal{H},D_{h}),   D_{h}=e^{h}De^{h}.
\end{equation}
The operator $D_{h}$ is still self-adjoint but the basic boundedness
condition
$[D,a]$  is bounded, for $\forall a \in \mathcal{A}$
will not necessarily hold, unless $h$ is in the center of $\mathcal{A}$.
In a twisted spectral triple \cite{Co7},
there is an automorphism $\sigma$ of the underlying noncommutative space $\mathcal{A}$
such that all twisted commutators $[D, a]_{\sigma}= D a- \sigma (a) D$ are bounded operators
($D$ is a selfadjoint operator playing the role of the Dirac operator).

In  \cite{Co8}, Connes and Chamseddine proved in the general framework of noncommutative geometry that the inner
fluctuations of the metric coming from the Morita equivalence ${\mathcal{A}%
}\sim{\mathcal{A}}$ generate perturbations of $D$ of the form $D\rightarrow
D^{\prime}=D+A$, where the $A$ plays the role of the gauge potentials and is a
self-adjoint element of the bimodule
\begin{equation}
\Omega_{D}^{1}=\,\Big\{\;\sum\,a_{j}\,[D,b_{j}]\;;\,a_{j},\,b_{j}\in{\mathcal{A}%
}\Big\}. \label{bim}%
\end{equation}
For any of the $a_{j}$ for $j>0$ is a scalar, it defines a functional on the universal $n$-forms $\Omega^{n}%
({\mathcal{A}})$ \cite{Co8} by the equality
\begin{equation}
\int_{\varphi}\,a_{0}\,da_{1}\cdots\,da_{n}=\;\varphi(a_{0},\,a_{1}%
,\cdots,a_{n}). \label{hoschtocyc}%
\end{equation}
The following functional is then a Hochschild
cocycle \cite{Co8} and is given as Dixmier trace of infinitesimals of order one,
\begin{equation}
\tau_{0}(a^{0},a^{1},a^{2},a^{3},a^{4})=\,{\int\!\!\!\!\!\!-}\,a^{0}%
\,[D,a^{1}]\,D^{-1}[D,a^{2}]\,D^{-1}[D,a^{3}]\,D^{-1}[D,a^{4}]\,D^{-1}.
\label{tau0}%
\end{equation}
The Dixmier trace \cite{JD} of an operator $T \in \mathcal{L}^{1,\infty}
(\mathcal{H}) $
  measures the {\it logarithmic divergence} of
its ordinary trace, and the Dixmier trace is invariant under
  perturbations by trace class operators. This is a very useful property and  makes the Dixmier trace
  a  flexible tool in computations.
In \cite{Co1,Co2}, Connes showed that the noncommutative residue on a compact manifold $M$ coincided with the Dixmier's trace on pseudo-differential operators of order $-{\rm {dim}}M$, and these non-normal traces can in fact be used to define
a process of noncommutative integration in noncommutative geometry, then a link between Wodzicki residue  and
the Dixmier's trace was given.
The Wodzicki residue in \cite{Wo,Wo1} gives the unique non-trivial trace
on the algebra of pseudodifferential operators,
And Connes \cite{Co2} claimed the Wodzicki residue (noncommutative residue) of the square of the inverse of the Dirac operator was proportioned to the Einstein-Hilbert action: $${\rm Wres}(D^{-2})=c_{0}\int_{M}s{\rm dvol}_{M},$$
where $c_0$ is a constant and $s$ is the scalar curvature, and $\upsilon_{n-1}=vol(S^{n-1})=\frac{2\pi^{m}}{\Gamma(m)}$ is the volume of the unit sphere $S^{n-1}$ in $\mathbb{R}^{n}$ and $n=2,m=4$.

Recently, using the Clifford representation of one-forms as $0$-order differential operators,
Dabrowski etc. \cite{DL} obtained the Einstein
tensor (or, more precisely, its contravariant version)   from functionals
over the dual bimodule of one-forms.  Dabrowski etc. \cite{DL} demonstrated that the noncommutative residue density recovered
the tensors $g$ and $G:=Ric-\frac{1}{2}R(g)g$ as certain bilinear functionals of vector fields on a manifold $M$, while their
dual tensors are recovered as a density of bilinear functionals of differential one-forms on $M$, which
recovered other important tensors in both the classical
setup as well as for the generalised or quantum geometries.
The notion of scalar curvature for spectral triples has also been formulated in
this manner \cite{Co4, Co5,CCS,Co6,LW,IL,B1,SCH}. For Riemannian manifold $M$ of even dimension $n = 2m$ equipped
with a metric tensor $g$ and the (scalar) Laplacian $\Delta$, a localised functional in $C^{\infty}(M)$ can be defined by
  \begin{equation}
  {\rm  Wres}(f\Delta^{-m+1})=\frac{n-2}{12}\upsilon_{n-1}\int_{M}fR(g){\rm vol}_{g},
\end{equation}
where $f\in C^{\infty}(M)$, $ R = R(g)$ is the scalar curvature, that is the $g$-trace $R=g^{jk}R_{jk }$ of
the Ricci tensor with components $R_{jk}$  in local coordinates,  $g^{jk}$ are the raised components
of the metric $g$.

Let $u$, $v$,$w$ with the components with
respect to local coordinates $u_{a}$,$v_{a}$ and $w_{a}$, respectively, be  differential forms represented in
such a way as endomorphisms (matrices) $c(u) $, $c(v) $ and $c(w) $ on the spinor bundle.
The torsion functional, as defined in \cite{DL2}, assigns to a triple of one-forms $(u,v,w)$:
\begin{equation}
	\mathscr{T}_D(u,v,w)= \mathrm{Wres}(c(u)c(v)c(w)D D^{-2m}), \qquad u,v,w\in \Omega^1_D.
\end{equation}
Finally, the scalar curvature functional is defined by $\mathscr{R}_D(f)=Wres(fD^{-2m+2})$ for $f\in \mathcal{A}$ \cite{Ka,KW}.
In \cite{WW2},  we give some new spectral functionals which is the extension of spectral functionals to the noncommutative realm with torsion,
 and we relate them to  the noncommutative residue for manifolds with boundary.
For a finitely summable regular spectral triple, we \cite{WW4} recover two forms, torsion of the linear connection and four forms by the noncommutative residue and perturbed de-Rham Hodge operators, and provide an explicit computation of generalized spectral forms associated with the perturbed de-Rham Hodge Dirac triple.
 Motivated by twisted spectral triple \cite{Co7}, the inner fluctuations of the spectral action \cite{Co8} and the spectral torsion \cite{IL} ,
the purpose of this paper is to generalize the results in \cite{DL,WW2,WW4,Wa3,WW6,WW3} and get
some new spectral torsion which is the extension of spectral torsion of Dirac operator to one form rescaled Dirac operator.
We will give the properties of the `perturbed' triple
\begin{equation}
(\mathcal{A},\mathcal{H},\widetilde{D} ),  \qquad \widetilde{D} =c(V)(D+\sqrt{-1}c(X))c(V),
\end{equation}
where the vector field $V$ is not zero at any point $x_{0}\in M$, and compute the spectral torsion $\mathscr{S}_{\widetilde{D}}$
\begin{equation}
\mathscr{S}_{\widetilde{D}}\big(c(u),c(v),c(w)\big)=\mathrm{Wres}\big(
\,a^{1}[\widetilde{D},a^{2}]\,a^{3}[\widetilde{D},a^{4}]\,a^{5}[\widetilde{D},a^{6}]\,\widetilde{D}  \,\widetilde{D}^{-2m}\big),
\label{tau0}%
\end{equation}
where $a^{j} \in C^{\infty}(M), j=1,...,6$.

The aim of this paper is to prove the following Theorem.
 \begin{thm}\label{mrthm}
 	Let $M$ be an $n=2m$ dimensional ($n\geq 3$) oriented compact spin Riemannian manifold,  with the trilinear Clifford multiplication by functional of differential one-forms $c(u),c(v),c(w),$ the spectral torsion $\mathscr{S}_{\widetilde{D}}$ for one-forms rescaled Dirac operator $\widetilde{D}=c(V)(D+\sqrt{-1}c(X))c(V)$ equals to
\begin{align}\label{sigma0}
&\mathscr{S}_{\widetilde{D}}\big(c(\widetilde{u})c(\widetilde{v})c(\widetilde{w})\big) \\
 =&Wres\Big( c(V)c(\widetilde{u})c(\widetilde{v})c(\widetilde{w})c(V)\big[c(V)(D+\sqrt{-1}c(X))c(V)\big]^{-2m+1}\Big) \nonumber\\
=&2^{m} \frac{2 \pi^{m}}{\Gamma (m)}\int_{M}\bigg\{
\frac{2m-3}{2}||V ||^{-4m+2}
\Big(g(\widetilde{u},\widetilde{v}) \widetilde{w} (||V ||^{2})
-g(\widetilde{u},\widetilde{w}) \widetilde{v} (||V ||^{2}))   +g(\widetilde{v},\widetilde{w}) \widetilde{u} (||V ||^{2})\Big)
\bigg\}{\rm d Vol}_M.\nonumber
\end{align}
 \end{thm}

\section{Preliminaries on the spectral torsion  for Dirac operator with one form rescaled}
In this section we fix notations and recall the previous work \cite{DL, WW4} that will play a fundamental role here.
We also give a review on the Lichnerowicz formula for one form rescaled Dirac operator  and the symbols
		of the  inverse of one form rescaled Dirac operator  and their relations.
\subsection{ The Lichnerowicz formula for Dirac operator with one form rescaled}
Let $M$ be a smooth compact Riemannian $n$-dimensional
 manifold without boundary and $V$ be a vector bundle on $M$. Recall that a differential operator $P$
 is of Laplace type if it has locally the form
 	\begin{align}
P=-(g^{ij}\partial_i\partial_j+A^i\partial_i+B),
	\end{align}
 where $\partial_i$ is a natural local frame on $TM$ and $g_{i,j}=g(\partial_i,\partial_j)$ and
 $(g^{ij})_{1\leq i,j\leq n}$ is the inverse matrix associated
 to the metric matrix $(g_{i,j})_{1\leq i,j\leq n}$ on $M$, and $A^i$ and $B$ are smooth sections of
 ${\rm End}(V)$ on $M$ (endomorphism). If $P$ is a Laplace type operator
 of the form (2.1), then (see \cite{WW4,WW6}) there is a unique connection
 $\nabla$ on $V$ and an unique endomorphism $E$ such that
 \begin{align}
P=-[g^{ij}(\nabla_{\partial_i}\nabla_{\partial_j}-\nabla_{\nabla^L_{\partial_i}{\partial_j}})+E] ,
	\end{align}
 where $\nabla^L$ denotes the Levi-civita connection on $M$.
  Now we let $M$ be a $n$-dimensional oriented spin manifold
with Riemannian metric $g$. The Dirac operator $D$ is
locally given as follows in terms of orthonormal frames $e_i,~1\leq
i\leq n$ and natural frames $\partial_i$ of $TM$, one has
	\begin{align}
D=\sum_{i,j}g^{ij}c(\partial_i)\nabla^S_{\partial_j}=\sum_{i}c(e_i)\nabla^S_{e_i},
	\end{align}
where $c(e_i)$ denotes the Clifford action which satisfies the relation
 $$c(e_i)c(e_j)+c(e_j)c(e_i)=-2\delta_i^j,$$ and
 	\begin{align}
\nabla^S_{\partial_i}=\partial_i+\sigma_i,~~\sigma_i=\frac{1}{4}\sum_{j,k}\left<\nabla^L_{\partial_i}e_j,e_k\right>c(e_j)c(e_k).
	\end{align}
In what follows, using the Einstein sum convention for repeated index summation:
	\begin{align}
\partial^j=g^{ij}\partial_i,~~\sigma^i=g^{ij}\sigma_j,~~\Gamma^k=g^{ij}\Gamma_{ij}^k.
	\end{align}
Recall the Lichnerowicz formula for the square of the Dirac operator, by (6a) in  \cite{Ka}, we have
	\begin{align}
D^2=-g^{ij}\partial_i\partial_j-2\sigma^j\partial_j+\Gamma^k\partial_k-g^{ij}[\partial_i(\sigma_j)+\sigma_i\sigma_j-\Gamma_{ij}^k\sigma_k]+\frac{1}{4}s,
	\end{align}
where $s$ is the scalar curvature.

\begin{prop}
The Lichnerowicz formula for Dirac operator with one form rescaled:
	\begin{align}\label{sigma0}
\widetilde{D}^2
=&-\|V\|^{2}c(V) \Big[ -g^{ij}c(V)\partial_i  \partial_j -2g^{ij} \partial_j(c(V))\partial_i- g^{ij} \partial_i[\partial_j(c(V))  ]                          -2\sigma^j[c(V)\partial_j +\partial_j(c(V))  ]   \\
       & +\Gamma^k[c(V)\partial_k +\partial_k(c(V))  ]-g^{ij}[\partial_i(\sigma_j)+\sigma_i\sigma_j-\Gamma_{ij}^k\sigma_k]c(V)+\frac{1}{4}sc(V) \Big]\nonumber \\
       &-\sqrt{-1}\|V\|^{2}c(V) \Big[\sum_{i,j}g^{ij}c(\partial_i) c(X)c(V)\partial_j+\sum_{i,j}g^{ij}c(\partial_i) \partial_j [c(X)c(V)]
              +\sum_{i,j}g^{ij}c(\partial_i) \sigma_j c(X)c(V) \Big] \nonumber \\
      &-\sqrt{-1}\|V\|^{2}c(V) c(X)\Big[ \sum_{i,j}g^{ij}c(\partial_i)  c(V)\partial_j+\sum_{i,j}g^{ij}c(\partial_i) \partial_j [ c(V)]
              +\sum_{i,j}g^{ij}c(\partial_i) \sigma_j  c(V)  \Big] \nonumber \\
     &-c(V) c({\rm d}(\|V\|^{2}))\Big[ \sum_{i,j}g^{ij}c(\partial_i)  c(V)\partial_j+\sum_{i,j}g^{ij}c(\partial_i) \partial_j [ c(V)]
              +\sum_{i,j}g^{ij}c(\partial_i) \sigma_j  c(V)  \Big] \nonumber \\
        & -\|V\|^{2}c(V) \|X\|^{2}c(V) -\sqrt{-1}c(V) c({\rm d}(\|V\|^{2}))c(X)c(V), \nonumber
\end{align}

where $s$ is the scalar curvature.
\end{prop}
\begin{proof}
Let $X$ be a
vector field on $M$ and we also denote the associated
Clifford action by $c(X)$. For  Dirac operator with one form rescaled $\widetilde{D}=c(V)(D+\sqrt{-1}c(X))c(V)$, then
	\begin{align}\label{sigma0}
\widetilde{D}^2=&c(V)(D+\sqrt{-1}c(X))c(V)c(V)(D+\sqrt{-1}c(X))c(V)\\
               =&-c(V)(D+\sqrt{-1}c(X))\|V\|^{2}(D+\sqrt{-1}c(X))c(V) \nonumber \\
               =&-c(V)\big[\|V\|^{2}D+c({\rm d}(\|V\|^{2}))+\sqrt{-1}\|V\|^{2}c(X)\big](D+\sqrt{-1}c(X))c(V)  \nonumber \\
                =&-c(V)\|V\|^{2}(D +\sqrt{-1} c(X))(D+\sqrt{-1}c(X))c(V) \nonumber \\
                  & -c(V) c({\rm d}(\|V\|^{2})) (D+\sqrt{-1}c(X))c(V) \nonumber \\
               =&-c(V)\|V\|^{2}\big[D^{2} +\sqrt{-1}D c(X)+\sqrt{-1} c(X)D-c^{2}(X)  \big]c(V) \nonumber \\
                  & -c(V) c({\rm d}(\|V\|^{2})) (D+\sqrt{-1}c(X))c(V) \nonumber \\
 =&-\|V\|^{2}c(V) D^{2}c(V) -\sqrt{-1}\|V\|^{2}c(V)D c(X)c(V)-\sqrt{-1}\|V\|^{2}c(V) c(X)Dc(V) \nonumber \\
                  & +\|V\|^{2}c(V)c^{2}(X)c(V)-c(V) c({\rm d}(\|V\|^{2})) D c(V)-\sqrt{-1}c(V) c({\rm d}(\|V\|^{2}))c(X)c(V). \nonumber
\end{align}
By the composition formula of pseudodifferential operators, we get
	\begin{align}\label{sigma0}
D^{2}c(V)=&-g^{ij}\partial_i\partial_jc(V)-2\sigma^j\partial_jc(V)+\Gamma^k\partial_kc(V)
            -g^{ij}[\partial_i(\sigma_j)+\sigma_i\sigma_j-\Gamma_{ij}^k\sigma_k]c(V)+\frac{1}{4}sc(V)\\
         =&-g^{ij}\partial_i [c(V)\partial_j +\partial_j(c(V))  ]
                  -2\sigma^j[c(V)\partial_j +\partial_j(c(V))  ]+\Gamma^k[c(V)\partial_k +\partial_k(c(V))  ]\nonumber \\
           & -g^{ij}[\partial_i(\sigma_j)+\sigma_i\sigma_j-\Gamma_{ij}^k\sigma_k]c(V)+\frac{1}{4}sc(V)\nonumber \\
 =&-g^{ij}c(V)\partial_i  \partial_j -g^{ij} \partial_i(c(V))\partial_j- g^{ij} \partial_i[\partial_j(c(V))  ]                          -g^{ij}\partial_i(c(V))\partial_j \nonumber \\
                  & -2\sigma^j[c(V)\partial_j +\partial_j(c(V))  ]+\Gamma^k[c(V)\partial_k +\partial_k(c(V))  ]\nonumber \\
           & -g^{ij}[\partial_i(\sigma_j)+\sigma_i\sigma_j-\Gamma_{ij}^k\sigma_k]c(V)+\frac{1}{4}sc(V)\nonumber \\
  =&-g^{ij}c(V)\partial_i  \partial_j -2g^{ij} \partial_j(c(V))\partial_i- g^{ij} \partial_i[\partial_j(c(V))  ]                          \nonumber \\
                  & -2\sigma^j[c(V)\partial_j +\partial_j(c(V))  ]+\Gamma^k[c(V)\partial_k +\partial_k(c(V))  ]\nonumber \\
       & -g^{ij}[\partial_i(\sigma_j)+\sigma_i\sigma_j-\Gamma_{ij}^k\sigma_k]c(V)+\frac{1}{4}sc(V).\nonumber
\end{align}
Also, straightforward computations yield
	\begin{align}\label{sigma0}
D c(X)c(V)=& \sum_{i,j}g^{ij}c(\partial_i)\nabla^S_{\partial_j}c(X)c(V)=  \sum_{i,j}g^{ij}c(\partial_i)(\partial_j+ \sigma_j)c(X)c(V)  \\
             =& \sum_{i,j}g^{ij}c(\partial_i) c(X)c(V)\partial_j+\sum_{i,j}g^{ij}c(\partial_i) \partial_j [c(X)c(V)]
              +\sum_{i,j}g^{ij}c(\partial_i) \sigma_j c(X)c(V), \nonumber
\end{align}
and
	\begin{align}\label{sigma0}
Dc(V)=& \sum_{i,j}g^{ij}c(\partial_i)\nabla^S_{\partial_j}c(X)= \sum_{i,j}g^{ij}c(\partial_i)(\partial_j+ \sigma_j) c(V)  \\
             =& \sum_{i,j}g^{ij}c(\partial_i)  c(V)\partial_j+\sum_{i,j}g^{ij}c(\partial_i) \partial_j [ c(V)]
              +\sum_{i,j}g^{ij}c(\partial_i) \sigma_j  c(V).\nonumber
\end{align}
Then we obtain
	\begin{align}\label{sigma0}
\widetilde{D}^2=&-\|V\|^{2}c(V) D^{2}c(V) -\sqrt{-1}\|V\|^{2}c(V)D c(X)c(V)-\sqrt{-1}\|V\|^{2}c(V) c(X)Dc(V) \\
                  & +\|V\|^{2}c(V)c^{2}(X)c(V)-c(V) c({\rm d}(\|V\|^{2})) D c(V)-\sqrt{-1}c(V) c({\rm d}(\|V\|^{2}))c(X)c(V) \nonumber \\
=&-\|V\|^{2}c(V) \Big[ -g^{ij}c(V)\partial_i  \partial_j -2g^{ij} \partial_j(c(V))\partial_i- g^{ij} \partial_i[\partial_j(c(V))  ]                          -2\sigma^j[c(V)\partial_j +\partial_j(c(V))  ]\nonumber \\
       & +\Gamma^k[c(V)\partial_k +\partial_k(c(V))  ]-g^{ij}[\partial_i(\sigma_j)+\sigma_i\sigma_j-\Gamma_{ij}^k\sigma_k]c(V)+\frac{1}{4}sc(V) \Big]\nonumber \\
       &-\sqrt{-1}\|V\|^{2}c(V) \Big[\sum_{i,j}g^{ij}c(\partial_i) c(X)c(V)\partial_j+\sum_{i,j}g^{ij}c(\partial_i) \partial_j [c(X)c(V)]
              +\sum_{i,j}g^{ij}c(\partial_i) \sigma_j c(X)c(V) \Big] \nonumber \\
      &-\sqrt{-1}\|V\|^{2}c(V) c(X)\Big[ \sum_{i,j}g^{ij}c(\partial_i)  c(V)\partial_j+\sum_{i,j}g^{ij}c(\partial_i) \partial_j [ c(V)]
              +\sum_{i,j}g^{ij}c(\partial_i) \sigma_j  c(V)  \Big] \nonumber \\
     &-c(V) c({\rm d}(\|V\|^{2}))\Big[ \sum_{i,j}g^{ij}c(\partial_i)  c(V)\partial_j+\sum_{i,j}g^{ij}c(\partial_i) \partial_j [ c(V)]
              +\sum_{i,j}g^{ij}c(\partial_i) \sigma_j  c(V)  \Big] \nonumber \\
        & +\|V\|^{2}c(V)c^{2}(X)c(V) -\sqrt{-1}c(V) c({\rm d}(\|V\|^{2}))c(X)c(V). \nonumber
\end{align}

\end{proof}

\subsection{ The algebra representation  of symbols for  the Dirac operator with one form rescaled}

By \cite{Wa3}, we have the Dirac operator
\begin{align}
D&=\sum^n_{i=1}c(e_i)\bigg(e_i-\frac{1}{4}\sum_{s,t}\omega_{s,t}
(e_i)c(e_s)c(e_t)\bigg).
\end{align}
The symbol expansion of a parametrix of $D$ is given,
\begin{align}\label{sigma0}
\sigma_{1}(D)=\sqrt{-1}c(\xi);~~~
		\sigma_{0}(D)=-\frac{1}{4}\sum_{i,s,t}\omega_{st}(e_i)c(e_i)c(e_s)c(e_t).
	\end{align}
For simplicity, we state the  computation of the leading coefficients of the symbols
of $\Delta^{-1}$ (we assume that the kernel of $\Delta $ is finite dimensional and
		can be neglected in the following) for a second-order differential operator
		$\Delta $, with the symbol expansion,
		%
		\begin{equation}
			\sigma (\Delta )(x,\xi )=p_{2}+p_{1}+p_{0}.
			\label{eqA.3}
		\end{equation}
		The inverse is a pseudodifferential operator $P^{-1}$, with a symbol of
		the form,
		%
		\begin{equation}
			\sigma (\Delta ^{-1})(x,\xi )=q_{-2}+q_{-3}+q_{-4}+\cdots,
			\label{eqA.4}
		\end{equation}
		where $r_{k}$ is homogeneous in $\xi $ of order $-k$. 	
By the composition formula of pseudo-differential operators in \cite{Wa3}, we have
\begin{align}
1&\sim \sigma(\Delta \circ \Delta ^{-1})\sim\sum_{\alpha}\frac{1}{\alpha!}\partial_{\xi}^{\alpha}[\sigma(\Delta )]
   D_{x}^{\alpha}[\sigma(\Delta ^{-1})]\\
&=(p_2+p_1+p_0)(q_{-2}+q_{-3}+q_{-4}+\cdots)\nonumber\\
&~~~+\sum_j(\partial_{\xi_j}p_2+\partial_{\xi_j}p_1+\partial_{\xi_j}p_0)(
D_{x_j}q_{-2}+D_{x_j}q_{-3}+D_{x_j}q_{-4}+\cdots)\nonumber\\
&~~~+\sum_{i,j}(\partial_{\xi_i}\partial_{\xi_j}p_2+\partial_{\xi_i}\partial_{\xi_j}p_1+\partial_{\xi_i}\partial_{\xi_j}p_0)(
D_{x_i}D_{x_j}q_{-2}+D_{x_i}D_{x_j}q_{-3}+D_{x_i}D_{x_j}q_{-4}+\cdots)\nonumber\\
&=p_2q_{-2}+(p_1q_{-2}+p_2q_{-3}+\sum_j\partial_{\xi_j}p_2D_{x_j}q_{-2})+(p_0q_{-2}+p_1q_{-3}+p_2q_{-4}\nonumber\\
&+\sum_j\partial_{\xi_j}p_2D_{x_j}q_{-3}+\sum_{i,j}\partial_{\xi_i}\partial_{\xi_j}p_2D_{x_i}D_{x_j}q_{-2})+\cdots.
\end{align}
Then we obtain
		\begin{equation}
			\begin{aligned}
				&p_{2} q_{-2}=1,
				\\
				&p_{1} q_{-2}+p_{2} q_{-3} -i
				\partial_{\xi_{j}}(p_{2})\partial _{x_{i}}(q_{-2})=0,
			\end{aligned}
			\label{eqA.5}
		\end{equation}
		which we solve recursively, obtaining
		%
		\begin{equation}
			\label{SymPinverse}
			\begin{aligned}
				&q_{-2}=p_{2}^{-1},
				\\
				& q_{-3}= -q_{-2} \big( p_{1}
				q_{-2} -i
				\partial_{\xi_{j}}(p_{2})\partial _{x_{i}}(q_{-2}) \big).
			\end{aligned}
			%
		\end{equation}
		
 Let $M$ be an even-dimensional compact Riemannian manifold $M$ with components
	of the metric $g$ given in chosen local coordinates by $g_{ab}$. The Laplace
	operator, which is densely defined on $L^{2}(M,vol_{g})$, is expressed
	as
	%
	\begin{equation}
		\Delta = - \frac{1}{\sqrt{\text{det}(g)}} \partial _{a} \bigl( \sqrt{
			\text{det}(g)} g^{ab} \partial _{b} \bigr),
		\label{eq2.1}
	\end{equation}
	where $g^{ab}$ is the inverse of the matrix $g_{ab}$. For details, see \cite{DL}, though we use a different convention here.
		The symbols of the differential operator $\Delta $ are:
	%
	\begin{equation}
		\label{LapSym}
		\begin{aligned}
			\sigma_{2}(\Delta)  =& g^{ab} \xi _{a} \xi _{b}, \qquad 	\sigma_{1}(\Delta)
			= & \frac{-i}{\sqrt{\text{det}(g)}} \partial _{a} \bigl( \sqrt{\text{det}(g)}
			g^{ab} \bigr) \xi _{b}, \qquad 	\sigma_{0}(\Delta)=0.
		\end{aligned}
		%
	\end{equation}
 The Taylor expansion of $g_{ab}$ denote by:
	%
	\begin{equation}
		g_{ab} = \delta _{ab} - \frac{1}{3} R_{acbd} x^{c} x^{d} + o({\mathbf{x^{2}}}),
		\label{gNorm}
	\end{equation}
	and
	\begin{equation}
		\sqrt{\text{det}(g)} = 1 - \frac{1}{6} \mathrm{Ric}_{ab} x^{a} x^{b} + o({
			\mathbf{x^{2}}}),
		\label{volNorm}
	\end{equation}
	where $R_{acbd}$ and $\mathrm{Ric}_{ab}$ are the components of the Riemann
	and Ricci tensor, respectively, at the point with ${\mathbf{x}}=0$ and we use
	the notation $o({\mathbf{x^{k}}})$ to denote that we expand a function up to
	the polynomial of order $k$ in the normal coordinates. The inverse metric
	is

	\begin{equation}
		g^{ab} = \delta _{ab} + \frac{1}{3} R_{acbd} x^{c} x^{d} + o({\mathbf{x^{2}}}),
		\label{gInverseNorm}
	\end{equation}
	where $\delta _{ab}$ and $\delta ^{ab}$ denote the Kronecker symbols.
	Then the symbols of the Laplace operator in normal coordinates
	are
	%
	\begin{equation}
		\label{LapSymNorm}
		\begin{aligned}
				\sigma_{2}(\Delta) =& \bigl( \delta _{ab} + \frac{1}{3} R_{acbd} x^{c} x^{d}
			\bigr) \xi _{a} \xi _{b} + o({\mathbf{x^{2}}}),
			\\
				\sigma_{1}(\Delta) = & \frac{2i}{3} \mathrm{Ric}_{ab} x^{a} \xi _{b} + o({
				\mathbf{x^{2}}}).
		\end{aligned}
		%
	\end{equation}
	Then, one has:
	%
	\begin{lem}\cite{DL}
		\label{lem2.1}
		In normal coordinates around a fixed point of the manifold $M$, the symbols
		of the inverse of the Laplace operator read
		%
		\begin{equation}
			\begin{aligned}
				&\sigma_{-2}(\Delta) = ||\xi ||^{-4} \bigl( \delta _{ab} - \frac{1}{3} R_{acbd}
				x^{c} x^{d} \bigr) \xi _{a} \xi _{b} + o(\mathbf{x^{2}}),
				\\
				& \sigma_{-3}(\Delta)  = -\frac{2i}{3} \mathrm{Ric}_{ab} x^{a} \xi _{b} ||
				\xi ||^{-4} + o(\mathbf{x}).
    \end{aligned}
			\label{normalny_laplasjan}
		\end{equation}
	\end{lem}
For Dirac operator with one form rescaled $\widetilde{D}=c(V)(D+\sqrt{-1}c(X))c(V)$,  we have
	\begin{align}\label{sigma0}
\widetilde{D}=&c(V)(D+\sqrt{-1}c(X))c(V) \\
                =&c(V)D c(V)+c(V) \sqrt{-1}c(X) c(V) \nonumber\\
       =&c(V)\sum_{i,j}g^{ij}c(\partial_i)\nabla^S_{\partial_j} c(V)+\sqrt{-1}c(V)c(X) c(V) \nonumber\\
       =&c(V)\sum_{i,j}g^{ij}c(\partial_i)(\partial_j+ \sigma_j ) c(V)+\sqrt{-1}c(V)c(X) c(V) \nonumber\\
        =&c(V)\sum_{i,j}g^{ij}c(\partial_i)\partial_j  c(V)
        +c(V)\sum_{i,j}g^{ij}c(\partial_i) \sigma_j   c(V)+\sqrt{-1}c(V)c(X) c(V) \nonumber\\
                =& \sum_{i,j}g^{ij}c(V)c(\partial_i)  c(V)\partial_j+\sum_{i,j}g^{ij}c(V)c(\partial_i) \partial_j [ c(V)]
              +\sum_{i,j}g^{ij}c(V)c(\partial_i) \sigma_j  c(V) \nonumber\\
              &+\sqrt{-1}c(V)c(X) c(V) .\nonumber
\end{align}
Then, one has:
\begin{lem}
		\label{lem2.1}
		In normal coordinates around a fixed point of the manifold $M$, the symbols
		of  the Dirac operator with one form rescaled read
		%
		\begin{equation}
			\begin{aligned}
				&\sigma_{0}(\widetilde{D}) = \sum_{i,j}g^{ij}c(V)c(\partial_i) \partial_j [ c(V)]
              +\sum_{i,j}g^{ij}c(V)c(\partial_i) \sigma_j  c(V)  +\sqrt{-1}c(V)c(X) c(V),
				\\
				& \sigma_{1}(\widetilde{D})= \sum_{i,j}g^{ij}\sqrt{-1}c(V)c(\partial_i)  c(V) \xi_j.
    \end{aligned}
			\label{normalny_laplasjan}
		\end{equation}
	\end{lem}
By Proposition 2.1, we get
\begin{lem}
		\label{lem2.1}
		In normal coordinates around a fixed point of the manifold $M$, the symbols
		of $\widetilde{D}^{2}$ read
		%
		\begin{equation}
			\begin{aligned}
				\sigma_{1}(\widetilde{D}^{2}) =&
2\sqrt{-1}|V|^{2}c(V)g^{ij} \partial_j(c(V))\xi_i  -2\sqrt{-1}\sigma^jc(V)\xi_j  +\sqrt{-1}\Gamma^kc(V)\xi_k\nonumber\\
       &+|V|^{2}c(V)  \sum_{i,j}g^{ij}c(\partial_i) c(X)c(V)\xi_j
            +|V|^{2}c(V) c(X) \sum_{i,j}g^{ij}c(\partial_i)  c(V)\xi_j\nonumber\\
     &-\sqrt{-1}c(V) c({\rm d}(|V|^{2}))  \sum_{i,j}g^{ij}c(\partial_i)  c(V)\xi_j,
				\\
				 \sigma_{2}(\widetilde{D}^{2}) =& ||V ||^{4} ||\xi ||^{2}.
    \end{aligned}
			\label{normalny_laplasjan}
		\end{equation}
	\end{lem}
By (2.20) and Lemma 2.4, we obtain
	  	\begin{lem}
		\label{lem2.1}
		In normal coordinates around a fixed point of the manifold $M$, the symbols
		of the inverse of $\widetilde{D}^{2}$ read
		%
		\begin{equation}
			\begin{aligned}
				 \sigma_{-3}(\widetilde{D}^{2})=&  -\|V\|^{-8} ||\xi ||^{-4}  \Big[ 2\sqrt{-1}\|V\|^{2}c(V)g^{ij} \partial_j(c(V))\xi_i  -2\sqrt{-1}\sigma^jc(V)\xi_j    +\sqrt{-1}\Gamma^kc(V)\xi_k   \nonumber\\
       &+\|V\|^{2}c(V)  \sum_{i,j}g^{ij}c(\partial_i) c(X)c(V)\xi_j
            +\|V\|^{2}c(V) c(X) \sum_{i,j}g^{ij}c(\partial_i)  c(V)\xi_j\nonumber\\
     &-\sqrt{-1}c(V) c({\rm d}(\|V\|^{2}))  \sum_{i,j}g^{ij}c(\partial_i)  c(V)\xi_j\Big]\nonumber\\
     &+2\sqrt{-1}\|\xi\|^{-4}\xi^\mu\partial_{x_\mu}(\|V\|^{-4})
     -2\sqrt{-1}\|V\|^{-4}\|\xi\|^{-6}\xi^\mu\xi_\alpha\xi_\beta\partial_{x_\mu}g^{\alpha\beta},
				\\
				 \sigma_{-2}(\widetilde{D}^{2}) =& \|V\|^{-4} ||\xi ||^{-2}.
    \end{aligned}
			\label{normalny_laplasjan}
		\end{equation}
	\end{lem}
 \section{The spectral torsion for  Dirac operator with one form rescaled}
This section is designed to get the spectral torsion for  Dirac operator with one form rescaled.

 \subsection{Torsion functional for  Dirac operator with one form rescaled}

Let $S^{*}M\subset T^{*}M$ denotes the co-sphere bundle on $M$
 and a pseudo-differential operator $P\in \Psi DO(E)$,
denote by $\sigma_{-n}^{P}$ the component of order
$-n$ of the complete symbol $\sigma^{P}= \sum_{i}\sigma_{i}^{P}$ of $P$
such that the equality
 \begin{equation}
{\rm  Wres}(P)= \int_{S^{*}M}\text{trace}(\sigma_{-n}^{P}(x,\xi)){\rm  d}x {\rm  d}\xi.
\end{equation}
In \cite{Wo,Wo1,Co2,Co3,Ka,KW}, it was shown that the noncommutative  residue ${\rm  Wres}(\Delta^{-n/2+1})$ of a generalized laplacian $\Delta$  on a complex vector bundle $E$ over a closed compact manifold $M$,
is the integral of the second coefficient of the heat kernel expansion of $\Delta$ up to a proportional factor.
In \cite{Co1}, the well-known Connes' trace theorem states the Dixmier trace of $-n$ order
pseudo-differential operator equals to its  noncommutative  residue up to a constant on
a closed $n-$dimensional manifold.
 Denote by $\Delta$ the Laplacian as above and $Tr_{\omega}$ the Dixmier trace,
then
 \begin{equation}
Tr_{\omega}((1+\Delta)^{-n/2})=\frac{1}{n}{\rm  Wres}((1+\Delta)^{-n/2})=\frac{1}{n}{\rm  dim}(E){\rm  Vol}(S^{n-1}){\rm  Vol}_{M}.
\end{equation}
\begin{defn}\cite{DL}
Let $c(u)=\sum_{r=1}^{n} u_{r}c(e_r), c(v)=\sum_{p=1}^{n} v_{p}c(e_{p}), c(w)=\sum_{q=1}^{n} w_{q}c(e_{q}),$ the trilinear Clifford multiplication by functional of differential one-forms $c(u), c(v)$, $ c(w)$
\begin{align}
\mathscr{T}_D(u,v,w)= \mathrm{Wres}(c(u)c(v)c(w)D D^{-2m})
 \end{align}
  is called torsion functional.
  \end{defn}
 Now for  the `perturbed' triple $(\mathcal{A},\mathcal{H},\widetilde{D} ),$
 one has
   \begin{lem}
By the trilinear Clifford multiplication by functional of differential one-forms $c(\widetilde{u}), c(\widetilde{v})$,  $ c(\widetilde{w})$, the spectral torsion $\mathscr{S}_{\widetilde{D}}$ for $\widetilde{D}=c(V)(D+\sqrt{-1}c(X))c(V)$ defined by
\begin{align}
 	\mathscr{S}_{\widetilde{D}}\bigg(c(\widetilde{u}),c(\widetilde{v}),c(\widetilde{w})\bigg)
 &=\mathrm{Wres}\bigg(c(V)c(\widetilde{u})c(\widetilde{v})c(\widetilde{w})c(V)\big[c(V)(D+\sqrt{-1}c(X))c(V)\big]^{-2m+1}\bigg).
 \end{align}
   \end{lem}
 \begin{proof}
 We note that $[\widetilde{D},a^{j}]=c(V)[D,a^{j}]c(V)=c(V)c(da^{j})c(V)$, $a^{j} \in C^{\infty}(M), j=1,...,6$, then
\begin{align}\label{sigma0}
&a^{1}[\widetilde{D},a^{2}]\,a^{3}[\widetilde{D},a^{4}]\,a^{5}[\widetilde{D},a^{6}]\,\widetilde{D} \\
=&a^{1}c(V)c(da^{2})c(V)a^{3}c(V)c(da^{4})c(V)a^{5}c(V)c(da^{6})c(V) c(V)(D+\sqrt{-1}c(X))c(V)\nonumber \\
=&c(V)a^{1}c(da^{2})a^{3}(-||V ||^{-2})c(da^{4})a^{5}(-||V ||^{-2})c(da^{6})c(V) c(V)(D+\sqrt{-1}c(X))c(V)\nonumber \\
=& c(V)c(\widetilde{u})c(\widetilde{v})c(\widetilde{w})c(V) \big[c(V)(D+\sqrt{-1}c(X))c(V)\big],\nonumber
\end{align}
where $c(\widetilde{u})= a^{1}c(da^{2})$,$c(\widetilde{v})=a^{3}(-||V ||^{-2})c(da^{4})$, $c(\widetilde{w}) =a^{5}(-||V ||^{-2})c(da^{6})$.
Then  the spectral torsion $\mathscr{S}_{\widetilde{D}}$  for $\widetilde{D}=c(V)(D+\sqrt{-1}c(X))c(V)$
\begin{align}\label{sigma0}
\mathscr{S}_{\widetilde{D}}\big(c(\widetilde{u})c(\widetilde{v})c(\widetilde{w})\big)=& Wres\big(
\,a^{1}[\widetilde{D},a^{2}]\,a^{3}[\widetilde{D},a^{4}]\,a^{5}[\widetilde{D},a^{6}]\,\widetilde{D}  \,\widetilde{D}^{-2m}\big)\\
=&Wres\Big( c(V)c(\widetilde{u})c(\widetilde{v})c(\widetilde{w})c(V)\big[c(V)(D+\sqrt{-1}c(X))c(V)\big]^{-2m+1}\Big).\nonumber
\end{align}
 \end{proof}
 Let  $n=2 m $, by (3.1), we need to compute
\begin{align}\label{abd}
\int_{M}\int_{\|\xi\|=1} \operatorname{tr}\left[\sigma_{-2 m}
\Big( c(V)c(\widetilde{u})c(\widetilde{v})c(\widetilde{w})c(V)\big[c(V)(D+\sqrt{-1}c(X))c(V)\big]^{-2m+1}\Big)
\right](x, \xi)\sigma(\xi)dx.
\end{align}
By (3.8) in \cite{WW6}, we get
\begin{align}\label{ABD}
& \sigma_{-2 m}\left(c(V)c(\widetilde{u})c(\widetilde{v})c(\widetilde{w})c(V)\Big(c(V)(D+\sqrt{-1}c(X))c(V)\Big)^{-2 m+1}\right)\nonumber\\
=&c(V)c(\widetilde{u})c(\widetilde{v})c(\widetilde{w})c(V)\sigma_{-2 m}\left(\Big(c(V)(D+\sqrt{-1}c(X))c(V)\Big)^{-2 m}\cdot\Big(c(V)(D+\sqrt{-1}c(X))c(V)\Big)\right) \nonumber\\
=&c(V)c(\widetilde{u})c(\widetilde{v})c(\widetilde{w})c(V)\bigg\{\sum_{|\alpha|=0}^{\infty} \frac{(-i)^{|\alpha|}}{\alpha!} \partial_{\xi}^{\alpha}\Big[\sigma\Big(\big(c(V)(D+\sqrt{-1}c(X))c(V)\big)^{-2 m}\Big)\Big] \nonumber\\
&\times \partial_{x}^{\alpha}\left[\sigma\left(\Big(c(V)(D+\sqrt{-1}c(X))c(V)\Big)\right)\right]\bigg\}_{-2 m} \nonumber\\
=&H_{1}+H_{2}+H_{3},\nonumber
\end{align}
where
 \begin{equation}
H_{1}=c(V)c(\widetilde{u})c(\widetilde{v})c(\widetilde{w})c(V)\sigma_{-2 m} \big(c(V)(D+\sqrt{-1}c(X))c(V)\big)^{-2 m} \sigma_{0}\big(c(V)(D+\sqrt{-1}c(X))c(V)\big),
\end{equation}
 \begin{equation}
H_{2}=c(V)c(\widetilde{u})c(\widetilde{v})c(\widetilde{w})c(V)\sigma_{-2 m-1} \big(c(V)(D+\sqrt{-1}c(X))c(V)\big)^{-2 m} \sigma_{1}\big(c(V)(D+\sqrt{-1}c(X))c(V)\big),
\end{equation}
\begin{align}
H_{3}=&c(V)c(\widetilde{u})c(\widetilde{v})c(\widetilde{w})c(V)(-\sqrt{-1}) \sum_{j=1}^{2m} \partial_{\xi_{j}}\Big[\sigma_{-2 m} \big(c(V)(D+\sqrt{-1}c(X))c(V)\big)^{-2 m} \Big] \\
&\times \partial_{x_{j}}\Big[\sigma_{1}\big(c(V)(D+\sqrt{-1}c(X))c(V)\big)\Big].\nonumber
\end{align}

  \subsection{The spectral torsion for  Dirac operator with one form rescaled}

 Based on the algorithm yielding the principal symbol of a product of pseudo-differential operators in terms of the principal symbols of the factors, we have
 \begin{align}\label{sigma0}
&\mathscr{S}_{\widetilde{D}}\big(c(\widetilde{u})c(\widetilde{v})c(\widetilde{w})\big) \\
=& Wres\big(
\,a^{1}[\widetilde{D},a^{2}]\,a^{3}[\widetilde{D},a^{4}]\,a^{5}[\widetilde{D},a^{6}]\,\widetilde{D}  \,\widetilde{D}^{-2m}\big)\nonumber\\
=&Wres\Big( c(V)c(\widetilde{u})c(\widetilde{v})c(\widetilde{w})c(V)\big[c(V)(D+\sqrt{-1}c(X))c(V)\big]^{-2m+1}\Big) \nonumber\\
=&\int_{M}\int_{\|\xi\|=1} \operatorname{tr}\left[\sigma_{-2 m}
\Big( c(V)c(\widetilde{u})c(\widetilde{v})c(\widetilde{w})c(V)\big[c(V)(D+\sqrt{-1}c(X))c(V)\big]^{-2m+1}\Big)
\right](x, \xi)\sigma(\xi)dx \nonumber\\
=&\int_{M}\int_{\|\xi\|=1} \operatorname{tr}\Big[H_{1} (x_{0} )+H_{2} (x_{0} )+H_{3} (x_{0} )
\Big](x, \xi)\sigma(\xi)dx. \nonumber
\end{align}
 Firstly, we review here technical tool of the computation, which are the integrals of polynomial functions over the unit spheres. By (32) in \cite{B1}, we define
\begin{align*}
I_{S_n}^{\gamma_1\cdot\cdot\cdot\gamma_{2\bar{n}+2}}=\int_{|x|=1}d^nxx^{\gamma_1}\cdot\cdot\cdot x^{\gamma_{2\bar{n}+2}},
\end{align*}
i.e. the monomial integrals over a unit sphere.
Then by Proposition A.2. in \cite{B1}, polynomial integrals over higher spheres in the $n$-dimesional case are given
\begin{align}
I_{S_n}^{\gamma_1\cdot\cdot\cdot\gamma_{2\bar{n}+2}}=\frac{1}{2\bar{n}+n}[\delta^{\gamma_1\gamma_2}I_{S_n}^{\gamma_3\cdot\cdot\cdot\gamma_{2\bar{n}+2}}+\cdot\cdot\cdot+\delta^{\gamma_1\gamma_{2\bar{n}+1}}I_{S_n}^{\gamma_2\cdot\cdot\cdot\gamma_{2\bar{n}+1}}],
\end{align}
where $S_n\equiv S^{n-1}$ in $\mathbb{R}^n$.
For $\bar{n}=0$, $I^0={\rm Vol}(S^{n-1})$=$\frac{2\pi^{\frac{n}{2}}}{\Gamma(\frac{n}{2})}$, and we obtain
\begin{align}
I_{S_n}^{\gamma_1\gamma_2}&=\frac{1}{n}{\rm Vol}(S^{n-1})\delta^{\gamma_1}_{\gamma_2}.
\end{align}

Next, we compute each term of (3.8),(3.9),(3.10) in turn.

 {\bf (I):} Explicit representation for the first item:$ \int_{\|\xi\|=1} \operatorname{tr} [H_{1} (x_{0} ) ](x, \xi)  \sigma(\xi) $

By lemma 2.3, lemma 2.4 and the composition formula of pseudo-differential operators, the following results are given.
\begin{align}
&\sigma_2^{(-m+1)}=||V ||^{-4m+4}\|\xi\|^{-2m+2};~~~(\sigma_2^{-1})^k=||V ||^{-4k}\|\xi\|^{-2k};\nonumber\\
&\partial_{\xi_\mu}\sigma_2^{(-m+k+1)}=2(-m+k+1)||V ||^{-4m+4k+4}\|\xi\|^{-2m+2k}\xi^\mu;\nonumber\\
&\partial_{x_\mu}\sigma_2^{-1}=\partial_{x_\mu}(||V ||^{-4})\|\xi\|^{-2}-||V ||^{-4}\|\xi\|^{-4}\xi_\alpha\xi_\beta\partial_{x_\mu}g^{\alpha\beta}.
\end{align}
In normal coordinates around a fixed point of the manifold $M$, $w_{s,t}(e_p)(x_0)=0$, then
\begin{align}\label{0-2m}
 &\sigma_{-2 m} \big(c(V)(D+\sqrt{-1}c(X))c(V)\big)^{-2 m} \sigma_{0}\big(c(V)(D+\sqrt{-1}c(X))c(V)\big)(x_{0})\\
 =&||V ||^{-4m}\|\xi\|^{-2m}\Big(\sum_{i,j}g^{ij}c(V)c(\partial_i) \partial_j [ c(V)]
              +\sum_{i,j}g^{ij}c(V)c(\partial_i) \sigma_j  c(V) \nonumber\\
              & +\sqrt{-1}c(V)c(X) c(V)\Big)(x_{0}) \nonumber\\
  =&||V ||^{-4m}\|\xi\|^{-2m}\Big(\sum_{i,j}g^{ij}c(V)c(\partial_i) \partial_j [ c(V)]
                +\sqrt{-1}c(V)c(X) c(V)\Big)(x_{0}) \nonumber\\
 =&\sum_{i,j}g^{ij}||V ||^{-4m}\|\xi\|^{-2m} c(V)c(\partial_i) \partial_j [ c(V)](x_{0})
                +\sqrt{-1}||V ||^{-4m}\|\xi\|^{-2m}c(V)c(X) c(V)\Big)(x_{0}).\nonumber
\end{align}
By  direct calculation, we obtain
\begin{align}\label{0-2m}
H_{1}(x_0)=& c(V)c(\widetilde{u})c(\widetilde{v})c(\widetilde{w})c(V)\sigma_{-2 m} \big(c(V)(D+\sqrt{-1}c(X))c(V)\big)^{-2 m} \\ &\times\sigma_{0}\big(c(V)(D+\sqrt{-1}c(X))c(V)\big)(x_0) \nonumber\\
 =&c(V)c(\widetilde{u})c(\widetilde{v})c(\widetilde{w})c(V)\sum_{i,j}g^{ij}||V ||^{-4m}\|\xi\|^{-2m} c(V)c(\partial_i)
 \partial_j [ c(V)](x_{0}) \nonumber\\
                &+\sqrt{-1}c(V)c(\widetilde{u})c(\widetilde{v})c(\widetilde{w})c(V)||V ||^{-4m}\|\xi\|^{-2m}c(V)c(X) c(V)\Big)(x_{0}).\nonumber
\end{align}
By the relation of the Clifford action, we  obtain
\begin{align}
	&{\rm tr}\Big(c(V)c(\widetilde{u})c(\widetilde{v})c(\widetilde{w})c(V) c(V)c(X) c(V)\Big)\\
=&||V ||^{4}{\rm tr}\Big( c(\widetilde{u})c(\widetilde{v})c(\widetilde{w}) c(X)  \Big)\nonumber\\
=&||V ||^{4}\Big(g(\widetilde{u},\widetilde{v})g(\widetilde{w},X)
-g(\widetilde{u},\widetilde{w})g(\widetilde{v},X)+g(\widetilde{u},X)g(\widetilde{v},\widetilde{w})\Big) {\rm tr}[id].\nonumber
\end{align}
On the other hand, by Lemma A.1 in \cite{JHW}, we have
\begin{lem} \cite{JHW}
For $X_1, X_2, X_3, X_4, X_5, X_6\in \Gamma(TM)$,
    \begin{align}
    	&{\bf (1)}\;{\rm tr}\Big[c(X_1)c( X_2)\Big]=-g(X_1,  X_2){\rm tr}[id];\label{trace2}\\
    	&{\bf (2)}\;{\rm tr}\Big[c(X_1)c( X_2)c( X_3)c( X_4)\Big]\label{trace4}\\
    	&=\Big[g(X_1,  X_4)g( X_2,  X_3)-g(X_1,  X_3)g( X_2,  X_4)+g(X_1,  X_2)g( X_3,  X_4)\Big]{\rm tr}[id];\nonumber\\
    	&{\bf (3)}\;{\rm tr}\Big[c(X_1)c( X_2)c( X_3)c( X_4)c( X_5)c( X_6)\Big]\label{trace6}\\
    	&=\Big[-g(X_1,  X_6)g ( X_2,  X_5)g ( X_3,  X_4)+g(X_1,  X_6)g ( X_2,  X_4)g ( X_3,  X_5)-g(X_1,  X_6)g ( X_2,  X_3)g ( X_4,  X_5)\nonumber\\
    	&+g(X_1, X_5)g ( X_2,  X_6)g ( X_3, X_4)-g(X_1, X_5)g ( X_2,  X_4)g ( X_3, X_6)+g(X_1, X_5)g ( X_2, X_3)g ( X_4, X_6)\nonumber\\
    	&-g(X_1, X_4)g ( X_2, X_6)g ( X_3, X_5)+g(X_1, X_4)g ( X_2, X_5)g ( X_3, X_6)-g(X_1, X_4)g ( X_2, X_3)g ( X_5, X_6)\nonumber\\
    	&+g(X_1, X_3)g ( X_2, X_6)g ( X_4, X_5)-g(X_1, X_3)g ( X_2, X_5)g ( X_4, X_6)+g(X_1, X_3)g ( X_2, X_4)g ( X_5, X_6)\nonumber\\
    	&-g(X_1, X_2)g ( X_3, X_6)g ( X_4, X_5)+g(X_1, X_2)g ( X_3, X_5)g ( X_4, X_6)-g(X_1, X_2)g ( X_3, X_4)g ( X_5, X_6)\Big]{\rm tr}[id].\nonumber
    \end{align}
\end{lem}
By (3) of Lemma 3.3, we obtain
\begin{align}
	&{\rm tr}\Big(c(V)c(\widetilde{u})c(\widetilde{v})c(\widetilde{w}) c(\partial_i)
  c(e_{\alpha}) \Big)\\
  &=\Big[-g(V,  e_{\alpha})g ( \widetilde{u},  \partial_i)g ( \widetilde{v}, \widetilde{w})
+g(V,  e_{\alpha})g ( \widetilde{u}, \widetilde{w})g ( \widetilde{v},  \partial_i)
-g(V,  e_{\alpha})g ( \widetilde{u}, \widetilde{v})g ( \widetilde{w},  \partial_i)\nonumber\\
    	&+g ( V,  \partial_i)g ( \widetilde{u},  e_{\alpha})g ( \widetilde{v}, \widetilde{w})
  -g (V,  \partial_i)g ( \widetilde{u}, \widetilde{w})g ( \widetilde{v},  e_{\alpha})
    +g (V,  \partial_i)g ( \widetilde{u}, \widetilde{v})g ( \widetilde{w},  e_{\alpha})\nonumber\\
&-g(V, \widetilde{w})g ( \widetilde{u}, e_{\alpha})g (  \widetilde{v}, \partial_i)
    +g(V, \widetilde{w})g ( \widetilde{u}, \partial_i)g ( \widetilde{v}, e_{\alpha})
    -g(V, \widetilde{w})g ( \widetilde{u}, \widetilde{v})g ( \partial_i, e_{\alpha})\nonumber\\
&+g(V, \widetilde{v})g ( \widetilde{u}, e_{\alpha})g (\widetilde{w}, \partial_i)
    -g(V, \widetilde{v})g (\widetilde{u}, \partial_i)g ( \widetilde{w}, e_{\alpha})
    +g(V, \widetilde{v})g ( \widetilde{u}, \widetilde{w})g ( \partial_i, e_{\alpha})\nonumber\\
&-g(V, \widetilde{u})g ( \widetilde{v},  e_{\alpha})g ( \widetilde{w}, \partial_i)
    +g(V, \widetilde{u})g ( \widetilde{v}, \partial_i)g ( \widetilde{w},  e_{\alpha})
    -g(V, \widetilde{u})g ( \widetilde{v}, \widetilde{w})g ( \partial_i, e_{\alpha})\Big]{\rm tr}[id].\nonumber
\end{align}

\begin{lem}The following equations hold:
     \begin{align}
{\bf (1)}&\sum_{i,j,\alpha}g^{ij} \partial_j [V_{\alpha}]\big[ -g(V,  e_{\alpha})g ( \widetilde{u},  \partial_i)g ( \widetilde{v}, \widetilde{w})  \big](x_{0})=-\frac{1}{2}\widetilde{u}(||V ||^{2} )g(\widetilde{v},\widetilde{w}) ; \nonumber\\
{\bf (2)}&\sum_{i,j,\alpha}g^{ij} \partial_j [V_{\alpha}]\big[g(V,  e_{\alpha})g ( \widetilde{u}, \widetilde{w})g ( \widetilde{v},  \partial_i)   \big](x_{0})=\frac{1}{2}\widetilde{v}(||V ||^{2} )g(\widetilde{u},\widetilde{w}) ; \nonumber\\
{\bf (3)}&\sum_{i,j,\alpha}g^{ij} \partial_j [V_{\alpha}]\big[-g(V,  e_{\alpha})g ( \widetilde{u}, \widetilde{v})g ( \widetilde{w},  \partial_i)\big](x_{0})=-\frac{1}{2}\widetilde{w}(||V ||^{2} )g(\widetilde{u},\widetilde{v}) ;\nonumber\\
{\bf (4)}&\sum_{i,j,\alpha}g^{ij} \partial_j [V_{\alpha}]\big[ g (V,  \partial_i)g ( \widetilde{u},  e_{\alpha})g ( \widetilde{v}, \widetilde{w}) \big](x_{0})=g(\widetilde{u},\nabla^{L}_{V}V )g(\widetilde{v},\widetilde{w});\nonumber\\
{\bf (5)}&\sum_{i,j,\alpha}g^{ij} \partial_j [V_{\alpha}]\big[  -g (V,  \partial_i)g ( \widetilde{u}, \widetilde{w})g ( \widetilde{v},  e_{\alpha}) \big](x_{0})=-g(\widetilde{v},\nabla^{L}_{V}V )g(\widetilde{u},\widetilde{w}); \nonumber\\
{\bf (6)}&\sum_{i,j,\alpha}g^{ij} \partial_j [V_{\alpha}]\big[g ( V,  \partial_i)g ( \widetilde{u}, \widetilde{v})g ( \widetilde{w},  e_{\alpha}) \big](x_{0})=-g(\widetilde{w},\nabla^{L}_{V}V )g(\widetilde{u},\widetilde{v}); \nonumber\\
{\bf (7)}&\sum_{i,j,\alpha}g^{ij} \partial_j [V_{\alpha}]\big[-g(V, \widetilde{w})g ( \widetilde{u}, e_{\alpha})g (  \widetilde{v}, \partial_i)   \big](x_{0})=-g(V,\widetilde{w})g(\widetilde{u},\nabla^{L}_{\widetilde{v}}V ); \nonumber\\
{\bf (8)}&\sum_{i,j,\alpha}g^{ij} \partial_j [V_{\alpha}]\big[g(V, \widetilde{w})g ( \widetilde{u}, \partial_i)g ( \widetilde{v}, e_{\alpha})   \big](x_{0})=g(V,\widetilde{w})g(\widetilde{v},\nabla^{L}_{\widetilde{u}}V ); \nonumber\\
{\bf (9)}&\sum_{i,j,\alpha}g^{ij} \partial_j [V_{\alpha}]\big[-g(V, \widetilde{w})g ( \widetilde{u}, \widetilde{v})g ( \partial_i, e_{\alpha})   \big](x_{0})=-\sum_{j}\partial_j(V_{j})g(V,\widetilde{w}) g(\widetilde{u},\widetilde{v});\nonumber\\
{\bf (10)}&\sum_{i,j,\alpha}g^{ij} \partial_j [V_{\alpha}]\big[g(V, \widetilde{v})g ( \widetilde{u}, e_{\alpha})g (\widetilde{w}, \partial_i)   \big](x_{0})=g(V,\widetilde{v})g(\widetilde{u},\nabla^{L}_{\widetilde{w}}V ); \nonumber\\
{\bf (11)}&\sum_{i,j,\alpha}g^{ij} \partial_j [V_{\alpha}]\big[ -g(V, \widetilde{v})g (\widetilde{u}, \partial_i)g ( \widetilde{w}, e_{\alpha})   \big](x_{0})=-g(V,\widetilde{v})g(\widetilde{w},\nabla^{L}_{\widetilde{u}}V ); \nonumber\\
{\bf (12)}&\sum_{i,j,\alpha}g^{ij} \partial_j [V_{\alpha}]\big[g(V, \widetilde{v})g ( \widetilde{u}, \widetilde{w})g ( \partial_i, e_{\alpha})   \big](x_{0})=-\sum_{j}\partial_j(V_{j})g(V,\widetilde{v}) g(\widetilde{u},\widetilde{w}); \nonumber\\
{\bf (13)}&\sum_{i,j,\alpha}g^{ij} \partial_j [V_{\alpha}]\big[-g(V, \widetilde{u})g ( \widetilde{v},  e_{\alpha})g ( \widetilde{w}, \partial_i) \big](x_{0})=g(V,\widetilde{u})g(\widetilde{v},\nabla^{L}_{\widetilde{w}}V ); \nonumber\\
{\bf (14)}&\sum_{i,j,\alpha}g^{ij} \partial_j [V_{\alpha}]\big[g(V, \widetilde{u})g ( \widetilde{v}, \partial_i)g ( \widetilde{w},  e_{\alpha})   \big](x_{0})=g(V,\widetilde{u})g(\widetilde{w},\nabla^{L}_{\widetilde{v}}V ); \nonumber\\
{\bf (15)}&\sum_{i,j,\alpha}g^{ij} \partial_j [V_{\alpha}]\big[ -g(V, \widetilde{u})g ( \widetilde{v}, \widetilde{w})g ( \partial_i, e_{\alpha}) \big](x_{0})=-\sum_{j}\partial_j(V_{j})g(V,\widetilde{u}) g(\widetilde{v},\widetilde{w}). \nonumber
    \end{align}
\end{lem}
\begin{proof}
{\bf (1)} In normal coordinates around a fixed point of the manifold $M$, we obtain
 \begin{align}
& \sum_{i,j,\alpha}g^{ij} \partial_j [V_{\alpha}]\big[ -g(V,  e_{\alpha})g ( \widetilde{u},  \partial_i)g ( \widetilde{v}, \widetilde{w})  \big](x_{0}) \nonumber\\
 =&  -\sum_{i,j,\alpha}g^{ij} \partial_j [V_{\alpha}]  g(V,  e_{\alpha})g ( \widetilde{u},  \partial_i)g ( \widetilde{v}, \widetilde{w}) (x_{0})                      \nonumber\\
  =&  -\sum_{ j,\alpha}  \partial_j [V_{\alpha}]  g(V,  e_{\alpha})g ( \widetilde{u},  \partial_j)g ( \widetilde{v}, \widetilde{w}) (x_{0})                      \nonumber\\
  =&  -\sum_{ j,\alpha}  g(V,  \partial_j [V_{\alpha}] e_{\alpha})g ( \widetilde{u},  \partial_j)g ( \widetilde{v}, \widetilde{w}) (x_{0})                      \nonumber\\
    =&  -\sum_{ j }  g(V,  \nabla ^{L} _{\partial_j} V)g ( \widetilde{u},  \partial_j)g ( \widetilde{v}, \widetilde{w}) (x_{0})                      \nonumber\\
      =&  -  g(V,  \nabla ^{L} _{\sum_{ j }g ( \widetilde{u},  \partial_j)\partial_j} V)g ( \widetilde{v}, \widetilde{w}) (x_{0})                      \nonumber\\
       =&  -  g(V,  \nabla ^{L} _{\widetilde{u}} V)g ( \widetilde{v}, \widetilde{w}) (x_{0})                      \nonumber\\
  =&-\frac{1}{2}\widetilde{u}(||V ||^{2} )g(\widetilde{v},\widetilde{w}). \nonumber
     \end{align}
 {\bf (4)} Also, straightforward computations yield
 \begin{align}
&\sum_{i,j,\alpha}g^{ij} \partial_j [V_{\alpha}]\big[ g (V,  \partial_i)g ( \widetilde{u},  e_{\alpha})g ( \widetilde{v}, \widetilde{w})\big](x_{0}) \nonumber\\
=&\sum_{ j,\alpha}  \partial_j [V_{\alpha}] g ( V,  \partial_j)g ( \widetilde{u},  e_{\alpha})g ( \widetilde{v}, \widetilde{w}) (x_{0})\nonumber\\
=&\sum_{ j,\alpha}  g (V,  \partial_j)g ( \widetilde{u},  \partial_j [V_{\alpha}]e_{\alpha})g ( \widetilde{v}, \widetilde{w}) (x_{0})\nonumber\\
=&\sum_{ j }  g ( V,  \partial_j)g ( \widetilde{u},  \nabla ^{L} _{\partial_j}V )g ( \widetilde{v}, \widetilde{w}) (x_{0})\nonumber\\
=&g ( \widetilde{u},  \nabla ^{L} _{\sum_{ j }  g ( V,  \partial_j)\partial_j}V )g ( \widetilde{v}, \widetilde{w}) (x_{0})\nonumber\\
=&g(\widetilde{u},\nabla^{L}_{V}V )g(\widetilde{v},\widetilde{w}).\nonumber
    \end{align}
     {\bf (9)} Similarly, we obtain
 \begin{align}
&\sum_{i,j,\alpha}g^{ij} \partial_j [V_{\alpha}]\big[-g(V, \widetilde{w})g ( \widetilde{u}, \widetilde{v})g ( \partial_i, e_{\alpha})   \big](x_{0})\nonumber\\
=&-\sum_{ j,\alpha}  \partial_j [V_{\alpha}]g(V, \widetilde{w})g ( \widetilde{u}, \widetilde{v})g ( \partial_j, e_{\alpha})(x_{0})   \nonumber\\
=&\sum_{ j,\alpha}  \partial_j [V_{\alpha}]g(V, \widetilde{w})g ( \widetilde{u}, \widetilde{v})\delta_{j} ^{\alpha  }(x_{0}) \nonumber\\
=&-\sum_{j}\partial_j(V_{j})g(V,\widetilde{w}) g(\widetilde{u},\widetilde{v}).\nonumber
    \end{align}
Other conclusions can be obtained through similar calculations
\end{proof}
Substituting into above calculation, we obtain
\begin{align}
	&{\rm tr}\Big(c(V)c(\widetilde{u})c(\widetilde{v})c(\widetilde{w})c(V)\sum_{i,j}g^{ij}c(V)c(\partial_i) \partial_j [ c(V)]\Big)(x_{0})\\
=&-\sum_{i,j}g^{ij}||V ||^{2}\partial_j [V_{\alpha}]{\rm tr}\Big(c(V)c(\widetilde{u})c(\widetilde{v})c(\widetilde{w}) c(\partial_i)
  c(e_{\alpha}) \Big)(x_{0}) \nonumber\\
=&- ||V ||^{2}\Big[   -\frac{1}{2}\widetilde{u}(||V ||^{2} )g(\widetilde{v},\widetilde{w})
 +\frac{1}{2}\widetilde{v}(||V ||^{2} )g(\widetilde{u},\widetilde{w})
 -\frac{1}{2}\widetilde{w}(||V ||^{2} )g(\widetilde{u},\widetilde{v})  \nonumber\\
&+ g(\widetilde{u},\nabla^{L}_{V}V )g(\widetilde{v},\widetilde{w})
 -g(\widetilde{v},\nabla^{L}_{V}V )g(\widetilde{u},\widetilde{w})
 +g(\widetilde{w},\nabla^{L}_{V}V )g(\widetilde{u},\widetilde{v})  \nonumber\\
 &-g(V,\widetilde{w})g(\widetilde{u},\nabla^{L}_{\widetilde{v}}V )
 +g(V,\widetilde{w})g(\widetilde{v},\nabla^{L}_{\widetilde{u}}V )
 -\sum_{j}\partial_j(V_{j})g(V,\widetilde{w}) g(\widetilde{u},\widetilde{v}) \nonumber\\
 &+g(V,\widetilde{v})g(\widetilde{u},\nabla^{L}_{\widetilde{w}}V )
 -g(V,\widetilde{v})g(\widetilde{w},\nabla^{L}_{\widetilde{u}}V )
 -\sum_{j}\partial_j(V_{j})g(V,\widetilde{v}) g(\widetilde{u},\widetilde{w})\nonumber\\
 &-g(V,\widetilde{u})g(\widetilde{v},\nabla^{L}_{\widetilde{w}}V )
+g(V,\widetilde{u})g(\widetilde{w},\nabla^{L}_{\widetilde{v}}V )
 -\sum_{j}\partial_j(V_{j})g(V,\widetilde{u}) g(\widetilde{v},\widetilde{w})\Big]{\rm tr}[id]. \nonumber
 \end{align}
By integrating formula we get
\begin{align}
& \int_{\|\xi\|=1} \operatorname{tr} [H_{1}  (x_{0} ) ](x, \xi)  \sigma(\xi)  \\
=&   \sqrt{-1}||V ||^{-4m+4} \Big(g(\widetilde{u},\widetilde{v})g(\widetilde{w},X)
-g(\widetilde{u},\widetilde{w})g(\widetilde{v},X)+g(\widetilde{u},X)g(\widetilde{v},\widetilde{w})\Big) {\rm tr}[id]{\rm Vol}(S^{2m-1})\nonumber\\
&-  ||V ||^{-4m+2}\Big[   -\frac{1}{2}\widetilde{u}(||V ||^{2} )g(\widetilde{v},\widetilde{w})
 +\frac{1}{2}\widetilde{v}(||V ||^{2} )g(\widetilde{u},\widetilde{w})
 -\frac{1}{2}\widetilde{w}(||V ||^{2} )g(\widetilde{u},\widetilde{v})  \nonumber\\
&+ g(\widetilde{u},\nabla^{L}_{V}V )g(\widetilde{v},\widetilde{w})
 -g(\widetilde{v},\nabla^{L}_{V}V )g(\widetilde{u},\widetilde{w})
 +g(\widetilde{w},\nabla^{L}_{V}V )g(\widetilde{u},\widetilde{v})  \nonumber\\
 &-g(V,\widetilde{w})g(\widetilde{u},\nabla^{L}_{\widetilde{v}}V )
 +g(V,\widetilde{w})g(\widetilde{v},\nabla^{L}_{\widetilde{u}}V )
 -\sum_{j}\partial_j(V_{j})g(V,\widetilde{w}) g(\widetilde{u},\widetilde{v}) \nonumber\\
 &+g(V,\widetilde{v})g(\widetilde{u},\nabla^{L}_{\widetilde{w}}V )
 -g(V,\widetilde{v})g(\widetilde{w},\nabla^{L}_{\widetilde{u}}V )
 -\sum_{j}\partial_j(V_{j})g(V,\widetilde{v}) g(\widetilde{u},\widetilde{w})\nonumber\\
 &-g(V,\widetilde{u})g(\widetilde{v},\nabla^{L}_{\widetilde{w}}V )
+g(V,\widetilde{u})g(\widetilde{w},\nabla^{L}_{\widetilde{v}}V )
 -\sum_{j}\partial_j(V_{j})g(V,\widetilde{u}) g(\widetilde{v},\widetilde{w})\Big]{\rm tr}[id]{\rm Vol}(S^{2m-1}). \nonumber
 \end{align}

 {\bf (II):} Explicit representation for the second item:$ \int_{\|\xi\|=1} \operatorname{tr} [H_{2}  (x_{0} ) ](x, \xi)  \sigma(\xi) $

Write $\sigma_2^{(-m+1)}:=\big[\sigma_{-2}\big(c(V)(D+\sqrt{-1}c(X))c(V)\big)^{-2}\big]^{m-1}$, then by (3.8) in \cite{JHW}, we have
\begin{align}\label{9000}
\sigma_{-2m-1}[\big(c(V)(D+\sqrt{-1}c(X))c(V)\big)^{-2m}]=&m\sigma_2^{(-m+1)}\sigma_{-3}\big[\big(c(V)(D+\sqrt{-1}c(X))c(V)\big)^{-2}\big]\\
   &-\sqrt{-1}\sum_{k=0}^{m-2}\sum_{\mu=1}^{2m}\partial_{\xi_\mu}\sigma_2^{(-m+k+1)}\partial_{x_\mu}\sigma_2^{-1}(\sigma_2^{-1})^k.\nonumber
\end{align}
By lemma 2.5 and the composition formula of pseudo-differential operators, the following results are given.
 	\begin{align}\label{sigma}
&\sigma_{-2m-1}[\big(c(V)(D+\sqrt{-1}c(X))c(V)\big)^{-2m}](x_{0})\\
=&    -m||V ||^{-4m-4}\|\xi\|^{-2m-2} \Big[ 2\sqrt{-1}|V|^{2}c(V)\sum_{i,j}g^{ij} \partial_j(c(V))\xi_i  -2\sigma^jc(V)\partial_j   +\Gamma^kc(V)\partial_k \nonumber\\
       &+|V|^{2}c(V)  \sum_{i,j}g^{ij}c(\partial_i) c(X)c(V)\xi_j
            +|V|^{2}c(V) c(X) \sum_{i,j}g^{ij}c(\partial_i)  c(V)\xi_j\nonumber\\
     &-\sqrt{-1}c(V) c({\rm d}(|V|^{2}))  \sum_{i,j}g^{ij}c(\partial_i)  c(V)\xi_j\Big](x_{0})\nonumber\\
     &+2\sqrt{-1}m ||V ||^{-4m+4}\|\xi\|^{-2m-2} \xi^\mu\partial_{x_\mu}(||V ||^{-4})(x_{0})\nonumber\\
     &-2\sqrt{-1}m||V ||^{-4m}\|\xi\|^{-2m-4}\xi^\mu\xi_\alpha\xi_\beta\partial_{x_\mu}g^{\alpha\beta}(x_{0})\nonumber\\
&-2\sqrt{-1}\sum_{k=0}^{m-2}\sum_{\mu=1}^{2m}||V ||^{-4m+4}(-m+k+1)\|\xi\|^{-2m-2}\xi^\mu\partial_{x_\mu}(||V ||^{-4})(x_{0})\nonumber\\
&+2\sqrt{-1}\sum_{k=0}^{m-2}\sum_{\mu=1}^{2m}
||V ||^{-4m}(-m+k+1)\|\xi\|^{-2m-4}\xi^\mu\xi_\alpha\xi_\beta\partial_{x_\mu}g^{\alpha\beta}(x_{0})\nonumber\\
=&    -m||V ||^{-4m-4}\|\xi\|^{-2m-2} \Big[ 2\sqrt{-1}|V|^{2}c(V)\sum_{i,j}g^{ij} \partial_j(c(V))\xi_i  \nonumber\\
       &+|V|^{2}c(V)  \sum_{i,j}g^{ij}c(\partial_i) c(X)c(V)\xi_j
            +|V|^{2}c(V) c(X) \sum_{i,j}g^{ij}c(\partial_i)  c(V)\xi_j\nonumber\\
     &-\sqrt{-1}mc(V) c({\rm d}(|V|^{2}))  \sum_{i,j}g^{ij}c(\partial_i)  c(V)\xi_j\Big](x_{0})\nonumber\\
     &+2\sqrt{-1} ||V ||^{-4m+4}\|\xi\|^{-2m-2} \xi^\mu\partial_{x_\mu}(||V ||^{-4})(x_{0})\nonumber\\
 &-2\sqrt{-1}\sum_{k=0}^{m-2}\sum_{\mu=1}^{2m}||V ||^{-4m+4}(-m+k+1)\|\xi\|^{-2m-2}\xi^\mu\partial_{x_\mu}(||V ||^{-4})(x_{0})\nonumber\\
 =&    -m||V ||^{-4m-2}\|\xi\|^{-2m-2}  2\sqrt{-1} c(V)\sum_{i,j}g^{ij} \partial_j(c(V))\xi_i  \nonumber\\
       & -m||V ||^{-4m-2}\|\xi\|^{-2m-2}  c(V) \Big[ \sum_{i,j}g^{ij}c(\partial_i) c(X)c(V)
             -  c(X) \sum_{i,j}g^{ij}c(\partial_i)  c(V)\Big]\xi_j\nonumber\\
     &+  m||V ||^{-4m-4}\|\xi\|^{-2m-2}  \sqrt{-1}c(V) c({\rm d}(|V|^{2}))  \sum_{i,j}g^{ij}c(\partial_i)  c(V)\xi_j (x_{0})\nonumber\\
     &+2\sqrt{-1} m||V ||^{-4m+4}\|\xi\|^{-2m-2} \xi^\mu\partial_{x_\mu}(||V ||^{-4})(x_{0})\nonumber\\
 &-2\sqrt{-1}\sum_{k=0}^{m-2}\sum_{\mu=1}^{2m}||V ||^{-4m+4}(-m+k+1)\|\xi\|^{-2m-2}\xi^\mu\partial_{x_\mu}(||V ||^{-4})(x_{0}).\nonumber
	\end{align}
In normal coordinates around a fixed point of the manifold $M$, we obtain
\begin{align}\label{0-2m}
 &\sigma_{-2 m-1} \big(c(V)(D+\sqrt{-1}c(X))c(V)\big)^{-2 m} \sigma_{1}\big(c(V)(D+\sqrt{-1}c(X))c(V)\big)(x_0)  \\
 =& 2m||V ||^{-4m-2}\|\xi\|^{-2m-2}  c(V)g^{ij} \partial_j(c(V))\xi_i \sum_{i,j}g^{ij} c(V)c(\partial_i)  c(V) \xi_j \nonumber\\
       & -\sqrt{-1}m||V ||^{-4m-2}\|\xi\|^{-2m-2}  c(V) \Big[ \sum_{i,j}g^{ij}c(\partial_i) c(X)c(V)
             -  c(X) \sum_{i,j}g^{ij}c(\partial_i)  c(V)\Big]\xi_ic(V)c(\partial_i)  c(V) \xi_j\nonumber\\
     &-  m||V ||^{-4m-4}\|\xi\|^{-2m-2}  c(V) c({\rm d}(|V|^{2}))  \sum_{i,j}g^{ij}c(\partial_i)  c(V)\xi_i\sum_{i,j}g^{ij} c(V)c(\partial_i)  c(V) \xi_j (x_{0})\nonumber\\
     &-2  m||V ||^{-4m+4}\|\xi\|^{-2m-2} \xi^\mu\partial_{x_\mu}(||V ||^{-4})\sum_{i,j}g^{ij} c(V)c(\partial_i)  c(V) \xi_j(x_{0})\nonumber\\
 &+2 \sum_{k=0}^{m-2}\sum_{\mu=1}^{2m}||V ||^{-4m+4}(-m+k+1)\|\xi\|^{-2m-2}\xi^\mu\partial_{x_\mu}(||V ||^{-4})\sum_{i,j}g^{ij} c(V)c(\partial_i)  c(V) \xi_j(x_{0}).\nonumber
\end{align}
Therefore
\begin{align}
H_{2}
=&c(V)c(\widetilde{u})c(\widetilde{v})c(\widetilde{w})c(V)\sigma_{-2 m-1} \big(c(V)(D+\sqrt{-1}c(X))c(V)\big)^{-2 m} \nonumber\\
    &\times \sigma_{1}\big(c(V)(D+\sqrt{-1}c(X))c(V)\big)  \nonumber\\
  =&L_{1}+L_{2}+L_{3}+L_{4}+L_{5},\nonumber
\end{align}
where
\begin{align}
L_{1}=& 2m||V ||^{-4m-2}\|\xi\|^{-2m-2}  c(V)c(\widetilde{u})c(\widetilde{v})c(\widetilde{w})c(V)c(V)\sum_{i,j}g^{ij} \partial_j(c(V))\xi_i g^{ij} c(V)c(\partial_i)  c(V) \xi_j; \nonumber\\
L_{2}=& -\sqrt{-1}m||V ||^{-4m-2}\|\xi\|^{-2m-2}c(V)c(\widetilde{u})c(\widetilde{v})c(\widetilde{w})c(V)  c(V) \nonumber\\
       &~~~\times\Big[ \sum_{i,j}g^{ij}c(\partial_i) c(X)c(V)
             -  c(X) \sum_{i,j}g^{ij}c(\partial_i)  c(V)\Big]\xi_ic(V)c(\partial_j)  c(V) \xi_j;\nonumber\\
L_{3}=&-  m||V ||^{-4m-4}\|\xi\|^{-2m-2} c(V)c(\widetilde{u})c(\widetilde{v})c(\widetilde{w})c(V) c(V) c({\rm d}(|V|^{2}))  \sum_{i,j}g^{ij}c(\partial_i)  c(V)\xi_i  c(V)c(\partial_j)  c(V) \xi_j;  \nonumber\\
L_{4}=&-2  m||V ||^{-4m+4}\|\xi\|^{-2m-2} \xi^\mu\partial_{x_\mu}(||V ||^{-4})\sum_{i,j}g^{ij}c(V)c(\widetilde{u})c(\widetilde{v})c(\widetilde{w})c(V) c(V)c(\partial_i)  c(V) \xi_j ;\nonumber\\
L_{5}=&2 \sum_{k=0}^{m-2}\sum_{\mu=1}^{2m}||V ||^{-4m+4}(-m+k+1)\|\xi\|^{-2m-2}c(V)c(\widetilde{u})c(\widetilde{v})c(\widetilde{w})c(V)\xi^\mu\nonumber\\
 &~~~~\times\partial_{x_\mu}(||V ||^{-4})\sum_{i,j}g^{ij} c(V)c(\partial_i)  c(V) \xi_j .\nonumber
\end{align}
By integrating formula we get
\begin{align}
& \int_{\|\xi\|=1} \operatorname{tr} [L_{1}  (x_{0} ) ](x, \xi)  \sigma(\xi)  \\
=&  \int_{\|\xi\|=1} \operatorname{tr} \Big[2m||V ||^{-4m-2}\|\xi\|^{-2m-2}  c(V)c(\widetilde{u})c(\widetilde{v})c(\widetilde{w})c(V)c(V)\sum_{i,j}g^{ij} \partial_j(c(V))\xi_i\nonumber\\
&~~~ \times \sum_{i,j}g^{ij} c(V)c(\partial_i)  c(V) \xi_j   (x_{0} ) \Big](x, \xi)  \sigma(\xi) \nonumber\\
=& 2m||V ||^{-4m+2}  \int_{\|\xi\|=1} \operatorname{tr} \Big[ c(\widetilde{u})c(\widetilde{v})c(\widetilde{w}) \sum_{i,j}g^{ij} \partial_j(c(V)) \xi_i   c(V)c(\partial_i)  c(V) \xi_j   (x_{0} ) \Big](x, \xi)  \sigma(\xi) \nonumber\\
=& 2m||V ||^{-4m+2}   \sum_{i,j}g^{ij}\operatorname{tr} \Big[ c(\widetilde{u})c(\widetilde{v})c(\widetilde{w}) \partial_j(c(V))    c(V)c(\partial_i)  c(V)     \Big] \frac{1}{2m} {\rm Vol}(S^{2m-1}) \nonumber\\
=&  ||V ||^{-4m+2}   \sum_{i,j,\alpha}g^{ij}\partial_j(c(V_{\alpha})) \operatorname{tr} \Big[ c(\widetilde{u})c(\widetilde{v})c(\widetilde{w}) c(e_{\alpha})c(V)c(\partial_i)  c(V)     \Big]  {\rm Vol}(S^{2m-1}) \nonumber\\
&=||V ||^{-4m+2} \sum_{i,j,\alpha}g^{ij}\partial_j(c(V_{\alpha})) \nonumber\\
&\times\Big[-g(\widetilde{u},   \partial_i)g ( \widetilde{v},  V)g ( \widetilde{w}, e_{\alpha})
+g(\widetilde{u},   \partial_i)g ( \widetilde{v}, e_{\alpha})g ( \widetilde{w},  V)
-g(\widetilde{u},   \partial_i)g ( \widetilde{v}, \widetilde{w})g ( e_{\alpha},  V)\nonumber\\
    	&+g ( \widetilde{u},  V)g ( \widetilde{v},  \partial_i))g ( \widetilde{w}, e_{\alpha})
  -g ( \widetilde{u},  V)g ( \widetilde{v}, e_{\alpha})g ( \widetilde{w},   \partial_i)
    +g ( \widetilde{u},  V)g ( \widetilde{v}, \widetilde{w})g (\partial_i,  e_{\alpha})\nonumber\\
&-g(\widetilde{u}, e_{\alpha})g ( \widetilde{v}, \partial_i)g (  \widetilde{w}, V)
    +g(\widetilde{u}, e_{\alpha})g ( \widetilde{v}, V)g ( \widetilde{w}, \partial_i)
    -g(\widetilde{u}, e_{\alpha})g ( \widetilde{v}, \widetilde{w})g ( V,\partial_i )\nonumber\\
&+g(\widetilde{u}, \widetilde{w})g ( \widetilde{v}, \partial_i)g (V, e_{\alpha})
    -g(\widetilde{u}, \widetilde{w})g (\widetilde{v}, V)g (\partial_i, e_{\alpha})
    +g(\widetilde{u}, \widetilde{w})g ( \widetilde{w}, e_{\alpha})g ( V,\partial_i )\nonumber\\
&-g(\widetilde{u}, \widetilde{v})g ( \widetilde{w},  \partial_i)g (V, e_{\alpha})
    +g(\widetilde{u}, \widetilde{v})g ( \widetilde{w}, V)g ( \partial_i,  e_{\alpha})
    -g(\widetilde{u}, \widetilde{v})g ( \widetilde{w}, e_{\alpha})g ( V,\partial_i)\Big]{\rm tr}[id]{\rm Vol}(S^{2m-1})\nonumber\\
=&||V ||^{-4m+2}\Big[   -\frac{1}{2}\widetilde{u}(||V ||^{2} )g(\widetilde{v},\widetilde{w})
 +\frac{1}{2}\widetilde{v}(||V ||^{2} )g(\widetilde{u},\widetilde{w})
 -\frac{1}{2}\widetilde{w}(||V ||^{2} )g(\widetilde{u},\widetilde{v})  \nonumber\\
&- g(\widetilde{u},\nabla^{L}_{V}V )g(\widetilde{v},\widetilde{w})
 +g(\widetilde{v},\nabla^{L}_{V}V )g(\widetilde{u},\widetilde{w})
 -g(\widetilde{w},\nabla^{L}_{V}V )g(\widetilde{u},\widetilde{v})  \nonumber\\
 &-g(V,\widetilde{w})g(\widetilde{u},\nabla^{L}_{\widetilde{v}}V )
 +g(V,\widetilde{w})g(\widetilde{v},\nabla^{L}_{\widetilde{u}}V )
 +\sum_{j}\partial_j(V_{j})g(V,\widetilde{w}) g(\widetilde{u},\widetilde{v}) \nonumber\\
 &+g(V,\widetilde{v})g(\widetilde{u},\nabla^{L}_{\widetilde{w}}V )
 -g(V,\widetilde{v})g(\widetilde{w},\nabla^{L}_{\widetilde{u}}V )
 -\sum_{j}\partial_j(V_{j})g(V,\widetilde{v}) g(\widetilde{u},\widetilde{w})\nonumber\\
 &-g(V,\widetilde{u})g(\widetilde{v},\nabla^{L}_{\widetilde{w}}V )
+g(V,\widetilde{u})g(\widetilde{w},\nabla^{L}_{\widetilde{v}}V )
 +\sum_{j}\partial_j(V_{j})g(V,\widetilde{u}) g(\widetilde{v},\widetilde{w})\Big]{\rm tr}[id]{\rm Vol}(S^{2m-1}). \nonumber
\end{align}
Repeated application of integrating formula yields that
\begin{align}
& \int_{\|\xi\|=1} \operatorname{tr} [L_{2}  (x_{0} ) ](x, \xi)  \sigma(\xi)  \\
=&  \int_{\|\xi\|=1} \operatorname{tr} \Big[ -\sqrt{-1}m||V ||^{-4m-2}\|\xi\|^{-2m-2}c(V)c(\widetilde{u})c(\widetilde{v})c(\widetilde{w})c(V)  c(V) \nonumber\\
       &~~~\times\Big( \sum_{i,j}g^{ij}c(\partial_i) c(X)c(V)
             -  c(X) \sum_{i,j}g^{ij}c(\partial_i)  c(V)\Big)\xi_ic(V)c(\partial_j)  c(V) \xi_j   (x_{0} ) \Big](x, \xi)  \sigma(\xi) \nonumber\\
=&  \sqrt{-1}m||V ||^{-4m+4} \sum_{i,j}g^{ij}\int_{\|\xi\|=1} \operatorname{tr} \Big[  c(\widetilde{u})c(\widetilde{v})c(\widetilde{w})   \Big(c(\partial_i) c(X)
             -  c(X)c(\partial_i)   \Big)\xi_i c(\partial_j)   \xi_j   \Big](x, \xi)  \sigma(\xi) \nonumber\\
=&  \sqrt{-1}m||V ||^{-4m+4} \sum_{i,j}g^{ij}  \operatorname{tr} \Big[  c(\widetilde{u})c(\widetilde{v})c(\widetilde{w})
\Big(c(\partial_i) c(X)-c(X)c(\partial_i)\Big) c(\partial_j)\Big]\frac{1}{2m} {\rm Vol}(S^{2m-1})\nonumber\\
=&  -2\sqrt{-1}m||V ||^{-4m+4} \sum_{i,j}g^{ij} g(X, \partial_i) \operatorname{tr} \Big[  c(\widetilde{u})c(\widetilde{v})c(\widetilde{w})
  c(\partial_j)\Big]\frac{1}{2m} {\rm Vol}(S^{2m-1})\nonumber\\
  =&  - \sqrt{-1} ||V ||^{-4m+4} \sum_{i,j}g^{ij} g(X, \partial_i) \Big(g(\widetilde{u},\widetilde{v})g(\widetilde{w},\partial_j)
-g(\widetilde{u},\widetilde{w})g(\widetilde{v},\partial_j)\nonumber\\
&~~~+g(\widetilde{v},\widetilde{w})g(\widetilde{u},\partial_j)\Big) {\rm tr}[id]  {\rm Vol}(S^{2m-1})\nonumber\\
 =&  - \sqrt{-1} ||V ||^{-4m+4}  \Big(g(\widetilde{u},\widetilde{v})g(\widetilde{w},X)
-g(\widetilde{u},\widetilde{w})g(\widetilde{v},X)
 +g(\widetilde{v},\widetilde{w})g(\widetilde{u},X)\Big) {\rm tr}[id]  {\rm Vol}(S^{2m-1});\nonumber
\end{align}
\begin{align}
& \int_{\|\xi\|=1} \operatorname{tr} [L_{3}  (x_{0} ) ](x, \xi)  \sigma(\xi)  \\
=&   \int_{\|\xi\|=1} \operatorname{tr} \Big[-  m||V ||^{-4m-4}\|\xi\|^{-2m-2} c(V)c(\widetilde{u})c(\widetilde{v})c(\widetilde{w})c(V) c(V) c({\rm d}(|V|^{2})) \nonumber\\
&~~~\times \sum_{i,j}g^{ij}c(\partial_i)  c(V)\xi_i  c(V)c(\partial_j)  c(V) \xi_j  (x_{0} ) \Big](x, \xi)  \sigma(\xi) \nonumber\\
=& m||V ||^{-4m+2}\int_{\|\xi\|=1} {\rm tr} \Big[  c(\widetilde{u})c(\widetilde{v})c(\widetilde{w}) c({\rm d}(|V|^{2}))
  \sum_{i,j}g^{ij}c(\partial_i)   \xi_i   c(\partial_j)  \xi_j  (x_{0} ) \Big](x, \xi)  \sigma(\xi) \nonumber\\
  =& -m||V ||^{-4m+2}  {\rm tr} \Big[  c(\widetilde{u})c(\widetilde{v})c(\widetilde{w}) c({\rm d}(|V|^{2}))
   \Big]   \frac{1}{2m} {\rm Vol}(S^{2m-1}) \nonumber\\
     =& \frac{- ||V ||^{-4m+2}}{2} \Big(g(\widetilde{u},\widetilde{v})g(\widetilde{w},{\rm d}(|V|^{2}))
-g(\widetilde{u},\widetilde{w})g(\widetilde{v},{\rm d}(|V|^{2}))
 +g(\widetilde{v},\widetilde{w})g(\widetilde{u},{\rm d}(|V|^{2}))\Big) {\rm tr}[id]    {\rm Vol}(S^{2m-1}); \nonumber
 \end{align}
and
\begin{align}
& \int_{\|\xi\|=1} {\rm tr} [L_{4}  (x_{0} ) ](x, \xi)  \sigma(\xi)  \\
=& \int_{\|\xi\|=1}{\rm tr}\Big[-2  m||V ||^{-4m+4}\|\xi\|^{-2m-2} \xi^\mu\partial_{x_\mu}(||V ||^{-4})\sum_{i,j}g^{ij}c(V)c(\widetilde{u})c(\widetilde{v})c(\widetilde{w})c(V) \nonumber\\
 & \times c(V)c(\partial_i)  c(V) \xi_j (x_{0} ) \Big](x, \xi)  \sigma(\xi)  \nonumber\\
 =& -2  m||V ||^{-4m+8}\int_{\|\xi\|=1} {\rm tr}\Big[\|\xi\|^{-2m-2} \xi^\mu\partial_{x_\mu}(||V ||^{-4})\sum_{i,j}g^{ij} c(\widetilde{u})c(\widetilde{v})c(\widetilde{w})
  c(\partial_i)    \xi_j (x_{0} ) \Big](x, \xi)  \sigma(\xi)  \nonumber\\
   =& -2  m||V ||^{-4m+8}\sum_{i,j}g^{ij} \partial_{x_\mu}(||V ||^{-4})\int_{\|\xi\|=1} {\rm tr} \Big[ c(\widetilde{u})c(\widetilde{v})c(\widetilde{w})
  c(\partial_i)    \xi^\mu\xi_j \Big](x_{0} )(x, \xi)  \sigma(\xi)  \nonumber\\
   =& -2  m||V ||^{-4m+8}\sum_{ j}  \partial_{x_j}(||V ||^{-4})  {\rm tr} \Big[ c(\widetilde{u})c(\widetilde{v})c(\widetilde{w})
  c(\partial_i)     \Big](x_{0} ) \frac{1}{2m} {\rm Vol}(S^{2m-1})  \nonumber\\
   =& - ||V ||^{-4m+8}\sum_{ j}  \partial_{x_j}(||V ||^{-4})\Big(g(\widetilde{u},\widetilde{v})g(\widetilde{w},\partial_j)
-g(\widetilde{u},\widetilde{w})g(\widetilde{v},\partial_j)\nonumber\\
 &+g(\widetilde{v},\widetilde{w})g(\widetilde{u},\partial_j)\Big) {\rm tr}[id]    {\rm Vol}(S^{2m-1})  \nonumber\\
  =& - ||V ||^{-4m+8} \Big(g(\widetilde{u},\widetilde{v}) \widetilde{w} (||V ||^{-4})
-g(\widetilde{u},\widetilde{w}) \widetilde{v} (||V ||^{-4})) +g(\widetilde{v},\widetilde{w}) \widetilde{u} (||V ||^{-4})\Big) {\rm tr}[id]    {\rm Vol}(S^{2m-1}).  \nonumber
 \end{align}
 Similarly, we obtain
\begin{align}
& \int_{\|\xi\|=1} {\rm tr}    [L_{5}  (x_{0} ) ](x, \xi)  \sigma(\xi)  \\
=&\int_{\|\xi\|=1} {\rm tr}     \Big[2 \sum_{k=0}^{m-2}\sum_{\mu=1}^{2m}||V ||^{-4m+4}(-m+k+1)\|\xi\|^{-2m-2}c(V)c(\widetilde{u})c(\widetilde{v})c(\widetilde{w})c(V)\xi^\mu\nonumber\\
 &~~~~\times\partial_{x_\mu}(||V ||^{-4})\sum_{i,j}g^{ij} c(V)c(\partial_i)  c(V) \xi_j (x_{0} ) \Big](x, \xi)  \sigma(\xi)  \nonumber\\
 =&2 \sum_{k=0}^{m-2}\sum_{\mu=1}^{2m}(-m+k+1)||V ||^{-4m+8}\partial_{x_\mu}(||V ||^{-4})\int_{\|\xi\|=1} {\rm tr}
 \Big[\|\xi\|^{-2m-2} c(\widetilde{u})c(\widetilde{v})c(\widetilde{w}) \xi^\mu\nonumber\\
 &~~~~\times\sum_{i,j}g^{ij} c(\partial_i)  \xi_j \Big](x_{0})(x, \xi)  \sigma(\xi)  \nonumber\\
  =&2 \sum_{k=0}^{m-2} (-m+k+1)||V ||^{-4m+8}\sum_{j}\partial_{x_j}(||V ||^{-4})  {\rm tr}
  \Big[  c(\widetilde{u})c(\widetilde{v})c(\widetilde{w})  \nonumber\\
 &~~~~\times c(\partial_j)   \Big](x_{0} ) \frac{1}{2m} {\rm Vol}(S^{2m-1})  \nonumber\\
  =&2 \sum_{k=0}^{m-2} (-m+k+1)||V ||^{-4m+8}\sum_{j}\partial_{x_j}(||V ||^{-4})  {\rm tr}
 \Big(g(\widetilde{u},\widetilde{v})g(\widetilde{w},\partial_j)
-g(\widetilde{u},\widetilde{w})g(\widetilde{v},\partial_j)\nonumber\\
 &+g(\widetilde{v},\widetilde{w})g(\widetilde{u},\partial_j)\Big)  \frac{1}{2m} {\rm Vol}(S^{2m-1})  \nonumber\\
  =&  \frac{1-m}{2}||V ||^{-4m+8}
\Big(g(\widetilde{u},\widetilde{v}) \widetilde{w} (||V ||^{-4})
-g(\widetilde{u},\widetilde{w}) \widetilde{v} (||V ||^{-4}))  \nonumber\\
&+g(\widetilde{v},\widetilde{w}) \widetilde{u} (||V ||^{-4})\Big) {\rm tr}[id]    {\rm Vol}(S^{2m-1}).  \nonumber
 \end{align}
Summing up (3.27)-(3.31) leads to the desired equality
\begin{align}
& \int_{\|\xi\|=1} \operatorname{tr} [H_{2}  (x_{0} ) ](x, \xi)  \sigma(\xi)  \\
=&||V ||^{-4m+2}\Big[   -\frac{1}{2}\widetilde{u}(||V ||^{2} )g(\widetilde{v},\widetilde{w})
 +\frac{1}{2}\widetilde{v}(||V ||^{2} )g(\widetilde{u},\widetilde{w})
 -\frac{1}{2}\widetilde{w}(||V ||^{2} )g(\widetilde{u},\widetilde{v})  \nonumber\\
&- g(\widetilde{u},\nabla^{L}_{V}V )g(\widetilde{v},\widetilde{w})
 +g(\widetilde{v},\nabla^{L}_{V}V )g(\widetilde{u},\widetilde{w})
 -g(\widetilde{w},\nabla^{L}_{V}V )g(\widetilde{u},\widetilde{v})  \nonumber\\
 &-g(V,\widetilde{w})g(\widetilde{u},\nabla^{L}_{\widetilde{v}}V )
 +g(V,\widetilde{w})g(\widetilde{v},\nabla^{L}_{\widetilde{u}}V )
 +\sum_{j}\partial_j(V_{j})g(V,\widetilde{w}) g(\widetilde{u},\widetilde{v}) \nonumber\\
 &+g(V,\widetilde{v})g(\widetilde{u},\nabla^{L}_{\widetilde{w}}V )
 -g(V,\widetilde{v})g(\widetilde{w},\nabla^{L}_{\widetilde{u}}V )
 -\sum_{j}\partial_j(V_{j})g(V,\widetilde{v}) g(\widetilde{u},\widetilde{w})\nonumber\\
 &-g(V,\widetilde{u})g(\widetilde{v},\nabla^{L}_{\widetilde{w}}V )
+g(V,\widetilde{u})g(\widetilde{w},\nabla^{L}_{\widetilde{v}}V )
 +\sum_{j}\partial_j(V_{j})g(V,\widetilde{u}) g(\widetilde{v},\widetilde{w})\Big]{\rm tr}[id]{\rm Vol}(S^{2m-1})\nonumber\\
 &- \sqrt{-1} ||V ||^{-4m+4}  \Big(g(\widetilde{u},\widetilde{v})g(\widetilde{w},X)
-g(\widetilde{u},\widetilde{w})g(\widetilde{v},X)
 +g(\widetilde{v},\widetilde{w})g(\widetilde{u},X)\Big) {\rm tr}[id]  {\rm Vol}(S^{2m-1}) \nonumber\\
&- \frac{||V ||^{-4m+2}}{2} \Big(g(\widetilde{u},\widetilde{v})g(\widetilde{w},{\rm d}(|V|^{2}))
-g(\widetilde{u},\widetilde{w})g(\widetilde{v},{\rm d}(|V|^{2}))
 +g(\widetilde{v},\widetilde{w})g(\widetilde{u},{\rm d}(|V|^{2}))\Big) {\rm tr}[id]    {\rm Vol}(S^{2m-1}) \nonumber\\
&-\frac{m+1}{2}||V ||^{-4m+8}
\Big(g(\widetilde{u},\widetilde{v}) \widetilde{w} (||V ||^{-4})
-g(\widetilde{u},\widetilde{w}) \widetilde{v} (||V ||^{-4}))   +g(\widetilde{v},\widetilde{w}) \widetilde{u} (||V ||^{-4})\Big) {\rm tr}[id]    {\rm Vol}(S^{2m-1}).  \nonumber
 \end{align}

 {\bf (III):} Explicit representation for the third item:$ \int_{\|\xi\|=1} \operatorname{tr} [H_{3}  (x_{0} ) ](x, \xi)  \sigma(\xi)$

Also, straightforward computations yield
\begin{align}\label{0-2m}
 & \partial_{\xi_{j}}\Big[\sigma_{-2 m} \big(c(V)(D+\sqrt{-1}c(X))c(V)\big)^{-2 m} \Big]
   \\
 =&\partial_{\xi_{j}}\Big[||V ||^{-4m}\|\xi\|^{-2m} \Big]=-2m||V ||^{-4m}\|\xi\|^{-2m-2}\xi^{\alpha},
\end{align}
and
\begin{align}\label{0-2m}
 & \partial_{x_{j}}\Big[\sigma_{1}\big(c(V)(D+\sqrt{-1}c(X))c(V)\big)\Big](x_{0}) \\
 =& \partial_{x_{j}}\Big[\sum_{i,j}g^{ij}\sqrt{-1} c(V)c(\partial_i)  c(V) \xi_j  \Big](x_{0}) \nonumber\\
 =& \sum_{i,j}g^{ij}\sqrt{-1}\partial_{x_{j}} (c(V))c(\partial_i)  c(V) \xi_j
+   \sum_{i,j}g^{ij}\sqrt{-1}c(V)c(\partial_i)  \partial_{x_{j}} (c(V))\xi_j .\nonumber
\end{align}
In normal coordinates around a fixed point of the manifold $M$, then
\begin{align}\label{0-2m}
 H_{3}(x_{0})=&c(V)c(\widetilde{u})c(\widetilde{v})c(\widetilde{w})c(V)
     (-\sqrt{-1}) \sum_{j=1}^{2m} \partial_{\xi_{j}}\Big[\sigma_{-2 m} \big(c(V)(D+\sqrt{-1}c(X))c(V)\big)^{-2 m} \Big]\\
  &\times \partial_{x_{j}}\Big[\sigma_{1}\big(c(V)(D+\sqrt{-1}c(X))c(V)\big)\Big] \nonumber\\
 =&c(V)c(\widetilde{u})c(\widetilde{v})c(\widetilde{w})c(V) (-\sqrt{-1})\sum_{j=1}^{2m} \big(-2m||V ||^{-4m}\|\xi\|^{-2m-2}\xi^{\alpha} \big)\nonumber\\
 &\times \big(\sum_{i,j}g^{ij}\sqrt{-1}\partial_{x_{j}} (c(V))c(\partial_i)  c(V) \xi_j
+ \sum_{i,j}g^{ij}\sqrt{-1}c(V)c(\partial_i)  \partial_{x_{j}} (c(V))\xi_j \big)\nonumber\\
=&K_{1}(x_{0})+K_{2} (x_{0}), \nonumber
\end{align}
where
\begin{align}\label{0-2m}
K_{1}(x_{0})=&-2m||V ||^{-4m}c(V)c(\widetilde{u})c(\widetilde{v})c(\widetilde{w})c(V)    \|\xi\|^{-2m-2}\xi^{\alpha}
  \sum_{i,j}g^{ij} \partial_{x_{j}} (c(V))c(\partial_i)  c(V) \xi_j;
\end{align}

\begin{align}\label{0-2m}
K_{2} (x_{0})=-2m||V ||^{-4m}\sum_{i,j}g^{ij}c(V)c(\widetilde{u})c(\widetilde{v})c(\widetilde{w})c(V)  \|\xi\|^{-2m-2}\xi^{\alpha}
 c(V)c(\partial_i)  \partial_{x_{j}} (c(V))\xi_j \big).
\end{align}
By integrating formula we get
\begin{align}
& \int_{\|\xi\|=1} \operatorname{tr} [K_{1}  (x_{0} ) ](x, \xi)  \sigma(\xi)  \\
=&  \int_{\|\xi\|=1} {\rm tr} \Big[-2m||V ||^{-4m}c(V)c(\widetilde{u})c(\widetilde{v})c(\widetilde{w})c(V)    \|\xi\|^{-2m-2}\xi^{\alpha}
  \sum_{i,j}g^{ij} \partial_{x_{j}} (c(V))c(\partial_i)  c(V) \xi_j \Big](x, \xi)  \sigma(\xi) \nonumber\\
  =&   2m||V ||^{-4m+2} \sum_{i,j,\beta}g^{ij}\partial_{x_{j}} (V_{\beta})\int_{\|\xi\|=1} {\rm tr} \Big[ c(\widetilde{u})c(\widetilde{v})c(\widetilde{w})c(V)    \xi^{\alpha} c(e_{\beta}) c(\partial_i)   \xi_j \Big](x, \xi)  \sigma(\xi) \nonumber\\
    =&   2m||V ||^{-4m+2} \sum_{ i,\beta} \partial_{x_{j}} (V_{\beta}) {\rm tr} \Big[ c(\widetilde{u})c(\widetilde{v})c(\widetilde{w})c(V)    c(e_{\beta}) c(\partial_i)\Big]  \frac{1}{2m} {\rm Vol}(S^{2m-1})     \nonumber\\
    =&    ||V ||^{-4m+2} \sum_{ i, \alpha} \partial_{x_{j}} (V_{\beta}) {\rm tr} \Big[ c(\widetilde{u})c(\widetilde{v})c(\widetilde{w})c(V)    c(e_{\alpha}) c(\partial_i)\Big]   {\rm Vol}(S^{2m-1})     \nonumber\\
      =&    ||V ||^{-4m+2}  \sum_{i, \alpha} \partial_j(c(V_{\alpha})) \nonumber\\
&\times\Big[-g(\widetilde{u},   \partial_i)g ( \widetilde{v}, e_{\alpha} )g ( \widetilde{w}, V)
+g(\widetilde{u},   \partial_i)g ( \widetilde{v},V)g ( \widetilde{w}, e_{\alpha})
-g(\widetilde{u},   \partial_i)g ( \widetilde{v}, \widetilde{w})g ( e_{\alpha},  V)\nonumber\\
    	&-g ( \widetilde{u},  V)g ( \widetilde{v},  \partial_i))g ( \widetilde{w}, e_{\alpha})
  +g ( \widetilde{u},  V)g ( \widetilde{v}, e_{\alpha})g ( \widetilde{w},   \partial_i)
   -g ( \widetilde{u},  V)g ( \widetilde{v}, \widetilde{w})g (\partial_i,  e_{\alpha})\nonumber\\
    &+g(\widetilde{u}, e_{\alpha})g ( \widetilde{v}, \partial_i)g (  \widetilde{w}, V)
    -g(\widetilde{u}, e_{\alpha})g ( \widetilde{v}, V)g ( \widetilde{w}, \partial_i)
   +g(\widetilde{u}, e_{\alpha})g ( \widetilde{v}, \widetilde{w})g ( V,\partial_i )\nonumber\\
    &+g(\widetilde{u}, \widetilde{w})g ( \widetilde{v}, \partial_i)g (V, e_{\alpha})
    -g(\widetilde{u}, \widetilde{w})g (\widetilde{v},e_{\alpha})g (\partial_i, V)
    +g(\widetilde{u}, \widetilde{w})g ( \widetilde{v},V)g ( e_{\alpha},\partial_i )\nonumber\\
    &-g(\widetilde{u}, \widetilde{v})g ( \widetilde{w},  \partial_i)g (V, e_{\alpha})
    +g(\widetilde{u}, \widetilde{v})g ( \widetilde{w},e_{\alpha})g ( \partial_i, V)
    -g(\widetilde{u}, \widetilde{v})g ( \widetilde{w},V )g ( e_{\alpha},\partial_i)\Big]{\rm tr}[id]{\rm Vol}(S^{2m-1})\nonumber\\
=&    ||V ||^{-4m+2} \Big[   -\frac{1}{2}\widetilde{u}(||V ||^{2} )g(\widetilde{v},\widetilde{w})
 +\frac{1}{2}\widetilde{v}(||V ||^{2} )g(\widetilde{u},\widetilde{w})
 -\frac{1}{2}\widetilde{w}(||V ||^{2} )g(\widetilde{u},\widetilde{v})  \nonumber\\
&+g(\widetilde{u},\nabla^{L}_{V}V )g(\widetilde{v},\widetilde{w})
 -g(\widetilde{v},\nabla^{L}_{V}V )g(\widetilde{u},\widetilde{w})
 +g(\widetilde{w},\nabla^{L}_{V}V )g(\widetilde{u},\widetilde{v})  \nonumber\\
 &+g(V,\widetilde{w})g(\widetilde{u},\nabla^{L}_{\widetilde{v}}V )
 -g(V,\widetilde{w})g(\widetilde{v},\nabla^{L}_{\widetilde{u}}V )
 -\sum_{j}\partial_j(V_{j})g(V,\widetilde{w}) g(\widetilde{u},\widetilde{v}) \nonumber\\
 &-g(V,\widetilde{v})g(\widetilde{u},\nabla^{L}_{\widetilde{w}}V )
 +g(V,\widetilde{v})g(\widetilde{w},\nabla^{L}_{\widetilde{u}}V )
 +\sum_{j}\partial_j(V_{j})g(V,\widetilde{v}) g(\widetilde{u},\widetilde{w})\nonumber\\
 &+g(V,\widetilde{u})g(\widetilde{v},\nabla^{L}_{\widetilde{w}}V )
-g(V,\widetilde{u})g(\widetilde{w},\nabla^{L}_{\widetilde{v}}V )
 -\sum_{j}\partial_j(V_{j})g(V,\widetilde{u}) g(\widetilde{v},\widetilde{w})\Big]{\rm tr}[id]{\rm Vol}(S^{2m-1}). \nonumber
\end{align}
Similarly, we obtain
\begin{align}
& \int_{\|\xi\|=1} \operatorname{tr} [K_{2}  (x_{0} ) ](x, \xi)  \sigma(\xi)  \\
=&  \int_{\|\xi\|=1} {\rm tr} \Big[-2m||V ||^{-4m}\sum_{i,j}g^{ij}c(V)c(\widetilde{u})c(\widetilde{v})c(\widetilde{w})c(V)  \|\xi\|^{-2m-2}\xi^{\alpha} c(V)c(\partial_i)  \partial_{x_{j}} (c(V))\xi_j \big)\Big](x, \xi)  \sigma(\xi) \nonumber\\
  =&   2m||V ||^{-4m+2} \sum_{i,j,\beta}g^{ij}\partial_{x_{j}} (V_{\beta})\int_{\|\xi\|=1} {\rm tr}
  \Big[ c(V) c(\widetilde{u})c(\widetilde{v})c(\widetilde{w})   \xi^{\alpha} c(\partial_i) c(e_{\beta}) \xi_j \Big](x, \xi)  \sigma(\xi) \nonumber\\
    =&   2m||V ||^{-4m+2} \sum_{ i,\beta} \partial_{x_{j}} (V_{\beta}) {\rm tr}
    \Big[ c(V) c(\widetilde{u})c(\widetilde{v})c(\widetilde{w})   c(\partial_i)c(e_{\beta})\Big]  \frac{1}{2m} {\rm Vol}(S^{2m-1})     \nonumber\\
    =&    ||V ||^{-4m+2} \sum_{ i, \alpha} \partial_{x_{j}} (V_{\beta}) {\rm tr}
     \Big[ c(V) c(\widetilde{u})c(\widetilde{v})c(\widetilde{w})  c(\partial_i) c(e_{\alpha})\Big]   {\rm Vol}(S^{2m-1})     \nonumber\\
      =&    ||V ||^{-4m+2}  \sum_{i, \alpha} \partial_j(c(V_{\alpha})) \nonumber\\
       & \times \Big[-g(V,  e_{\alpha})g ( \widetilde{u},  \partial_i)g ( \widetilde{v}, \widetilde{w})
+g(V,  e_{\alpha})g ( \widetilde{u}, \widetilde{w})g ( \widetilde{v},  \partial_i)
-g(V,  e_{\alpha})g ( \widetilde{u}, \widetilde{v})g ( \widetilde{w},  \partial_i)\nonumber\\
    	&+g ( V,  \partial_i)g ( \widetilde{u},  e_{\alpha})g ( \widetilde{v}, \widetilde{w})
  -g ( V,  \partial_i)g ( \widetilde{u}, \widetilde{w})g ( \widetilde{v},  e_{\alpha})
    +g (V,  \partial_i)g ( \widetilde{u}, \widetilde{v})g ( \widetilde{w},  e_{\alpha})\nonumber\\
&-g(V, \widetilde{w})g ( \widetilde{u}, e_{\alpha})g (  \widetilde{v}, \partial_i)
    +g(V, \widetilde{w})g ( \widetilde{u}, \partial_i)g ( \widetilde{v}, e_{\alpha})
    -g(V, \widetilde{w})g ( \widetilde{u}, \widetilde{v})g ( \partial_i, e_{\alpha})\nonumber\\
&+g(V, \widetilde{v})g ( \widetilde{u}, e_{\alpha})g (\widetilde{w}, \partial_i)
    -g(V, \widetilde{v})g (\widetilde{u}, \partial_i)g ( \widetilde{w}, e_{\alpha})
    +g(V, \widetilde{v})g ( \widetilde{u}, \widetilde{w})g ( \partial_i, e_{\alpha})\nonumber\\
&-g(V, \widetilde{u})g ( \widetilde{v},  e_{\alpha})g ( \widetilde{w}, \partial_i)
    +g(V, \widetilde{u})g ( \widetilde{v}, \partial_i)g ( \widetilde{w},  e_{\alpha})
    -g(V, \widetilde{u})g ( \widetilde{v}, \widetilde{w})g ( \partial_i, e_{\alpha})\Big]{\rm tr}[id]{\rm Vol}(S^{2m-1})\nonumber\\
   =&    ||V ||^{-4m+2}
\Big[   -\frac{1}{2}\widetilde{u}(||V ||^{2} )g(\widetilde{v},\widetilde{w})
 +\frac{1}{2}\widetilde{v}(||V ||^{2} )g(\widetilde{u},\widetilde{w})
 -\frac{1}{2}\widetilde{w}(||V ||^{2} )g(\widetilde{u},\widetilde{v})  \nonumber\\
&+ g(\widetilde{u},\nabla^{L}_{V}V )g(\widetilde{v},\widetilde{w})
 -g(\widetilde{v},\nabla^{L}_{V}V )g(\widetilde{u},\widetilde{w})
 +g(\widetilde{w},\nabla^{L}_{V}V )g(\widetilde{u},\widetilde{v})  \nonumber\\
 &-g(V,\widetilde{w})g(\widetilde{u},\nabla^{L}_{\widetilde{v}}V )
 +g(V,\widetilde{w})g(\widetilde{v},\nabla^{L}_{\widetilde{u}}V )
 -\sum_{j}\partial_j(V_{j})g(V,\widetilde{w}) g(\widetilde{u},\widetilde{v}) \nonumber\\
 &+g(V,\widetilde{v})g(\widetilde{u},\nabla^{L}_{\widetilde{w}}V )
 -g(V,\widetilde{v})g(\widetilde{w},\nabla^{L}_{\widetilde{u}}V )
 -\sum_{j}\partial_j(V_{j})g(V,\widetilde{v}) g(\widetilde{u},\widetilde{w})\nonumber\\
 &-g(V,\widetilde{u})g(\widetilde{v},\nabla^{L}_{\widetilde{w}}V )
+g(V,\widetilde{u})g(\widetilde{w},\nabla^{L}_{\widetilde{v}}V )
 -\sum_{j}\partial_j(V_{j})g(V,\widetilde{u}) g(\widetilde{v},\widetilde{w})\Big]{\rm tr}[id]{\rm Vol}(S^{2m-1}). \nonumber
 \end{align}
 From (3.39) and (3.40) we obtain
 \begin{align}
& \int_{\|\xi\|=1} \operatorname{tr} [H_{3}  (x_{0} ) ](x, \xi)  \sigma(\xi)  \\
=&    ||V ||^{-4m+2}
\Big[   - \widetilde{u}(||V ||^{2} )g(\widetilde{v},\widetilde{w})
 + \widetilde{v}(||V ||^{2} )g(\widetilde{u},\widetilde{w})
 - \widetilde{w}(||V ||^{2} )g(\widetilde{u},\widetilde{v})  \nonumber\\
&+ 2g(\widetilde{u},\nabla^{L}_{V}V )g(\widetilde{v},\widetilde{w})
 -2g(\widetilde{v},\nabla^{L}_{V}V )g(\widetilde{u},\widetilde{w})
 +2g(\widetilde{w},\nabla^{L}_{V}V )g(\widetilde{u},\widetilde{v})  \nonumber\\
 &  -2\sum_{j}\partial_j(V_{j})g(V,\widetilde{w}) g(\widetilde{u},\widetilde{v})
  -2\sum_{j}\partial_j(V_{j})g(V,\widetilde{u}) g(\widetilde{v},\widetilde{w})\Big]{\rm tr}[id]{\rm Vol}(S^{2m-1}). \nonumber
 \end{align}

Summing up
 {\bf (I)},
 {\bf (II)}, and
 {\bf (III)} leads to the desired equality

 \begin{align}\label{sigma0}
&\mathscr{S}_{\widetilde{D}}\big(c(\widetilde{u})c(\widetilde{v})c(\widetilde{w})\big) \\
 =&Wres\Big( c(V)c(\widetilde{u})c(\widetilde{v})c(\widetilde{w})c(V)\big[c(V)(D+\sqrt{-1}c(X))c(V)\big]^{-2m+1}\Big) \nonumber\\
 =&\int_{M}\int_{\|\xi\|=1} \operatorname{tr}\Big[H_{1} (x_{0} )+H_{2} (x_{0} )+H_{3} (x_{0} )
\Big](x, \xi)\sigma(\xi)dx \nonumber\\
=&2^{m} \frac{2 \pi^{m}}{\Gamma (m)}\int_{M}\bigg\{
\frac{2m-3}{2}||V ||^{-4m+2}
\Big(g(\widetilde{u},\widetilde{v}) \widetilde{w} (||V ||^{2})
-g(\widetilde{u},\widetilde{w}) \widetilde{v} (||V ||^{2}))   +g(\widetilde{v},\widetilde{w}) \widetilde{u} (||V ||^{2})\Big)
\bigg\}{\rm d Vol}_M.\nonumber
\end{align}
The proof of Theorem 1.1 is complete.

\section{Conclusions and outlook }
The basic paradigm of noncommutative geometry is that of a spectral triple  $(\mathcal{A},\mathcal{H},D)$, where the algebra  $ \mathcal{A}$ encodes the space and the operator  $ D $ encodes the metric.
The main point of this note, besides giving
simple natural examples of the general spectral torsion and developing  the
 rescaled spectral triple theory.
Based on the  noncommutative  residue methods, we represented the symbols of the inverse of the one form rescaled Dirac operator and the trace over endomorphisms of the bundle $E$ at given point of $M$. This approach connects spectral theory and noncommutative geometry, offering a deep insight into the geometric structure of both classical and noncommutative spaces.

\section*{ AUTHOR DECLARATIONS}
Conflict of Interest:

The authors have no conflicts to disclose.
\section*{ Author Contributions}
Jian Wang: Investigation (equal); Writing- original draft (equal). Yong Wang: Investigation (equal); Writing-original draft (equal).
\section*{DATA AVAILABILITY}
Data sharing is not applicable to this article as no new data were created or analyzed in this study.

\section*{ Acknowledgements}
The first author was supported by NSFC. 11501414. The second author was supported by NSFC. 11771070.
The authors also thank the referee for his (or her) careful reading and
helpful comments.

\end{document}